\documentclass[12pt,reqno]{amsart}
\usepackage{latexsym,verbatim,xypic}
\def\suffix{ps}
\newcount\system
\global\system=3   % for unix(dvips)
\def\ifundefined#1{\expandafter\ifx\csname#1\endcsname\relax}
\ifundefined{figdir}\def\figdir{}\fi
\newcount\firstline
\newdimen\pswidth  \newdimen\xleft
\newdimen\psheight \newdimen\ytop \newdimen\ybot
\newcount\justx \newcount\justy
\global\justx=0 \global\justy=0
\newdimen\vpos \newtoks\labeL 
\newread\labeLfile \newdimen\xcoord \newdimen\ycoord
\newif\ifdoit 
\newbox\labox
\newdimen\xdvikwid 
\newdimen\xdvikht
\newdimen\pspoints
\newdimen\rwi
\newcount\temp
\def\readdim#1{\global\read\labeLfile to \temp
\global #1=\temp pt}
\def\figcrop#1{\par%  #1=filename
\openin\labeLfile=\figdir#1.lbl                                              
\global\read\labeLfile to\firstline\message{#1}               
\global\read\labeLfile to\temp%read overall dimensions                                     
\readdim{\ybot}
\readdim{\xleft}%               read upper left point
\readdim{\ytop}
\global\read\labeLfile to\justx%ignore
\global\read\labeLfile to\justy%ignore
\global\read\labeLfile to\labeL%ignore
\readdim{\pswidth}%            read lower right point
\global\advance\pswidth by -\xleft
\readdim{\psheight}
\global\advance\ybot by -\psheight
\global\advance\psheight by -\ytop
\global\read\labeLfile to\justx%ignore
\global\read\labeLfile to\justy%ignore
\global\read\labeLfile to\labeL%ignore                                    
\vbox to\psheight{\vfill
\ifnum\system=1
\fi %textures
\ifnum\system=2
\fi %msdos
\ifnum\system=3
                                                 \fi         %%unix:dvips
\ifnum\system=4
\fi         %%unix:dvips,scaled
\ifnum\system=1
\hbox to \pswidth{\kern-\xleft\special{postscriptfile \figdir#1.\suffix }\hfil}\fi
\ifnum\system=2
\hbox to \pswidth{\kern-\xleft\special{ps: plotfile \figdir#1.\suffix }\hfil}\fi
\ifnum\system=3
\hbox to \pswidth{\kern-\xleft\includegraphics{\figdir#1.\suffix}\hfil}\fi
\ifnum\system=4
\hbox to \pswidth{\kern-\xleft\includegraphics{\figdir#1.\suffix}\hfil}\fi
\ifnum\system=5
\hbox to \pswidth{\kern-\xleft\includegraphics{\figdir#1.\suffix}\hfil}\fi %orphee
\ifnum\system=6
   \xdvikwid=\pswidth
   \xdvikht=\psheight
   {\global\divide\xdvikwid by \pspoints}
   {\global\divide\xdvikht by \pspoints}
   \rwi=\xdvikwid
    {\global\multiply\rwi by 10}
\hbox to \pswidth{\kern-\xleft\includegraphics{\figdir#1.\suffix\space}\hfil}\fi                   %xdvik
\vskip -\baselineskip
\vskip -\ybot 
\vskip-\psheight %                                     
\hbox to\pswidth  {\hss}%                                            
\parindent=0pt\offinterlineskip                                       
\vpos=0 pt%                                                              
\loop\readdim{\xcoord}                                 
\ifdim \xcoord < -999pt \doitfalse\else\doittrue\fi                        
\ifdoit \advance \xcoord by -\xleft
\readdim{\ycoord}
\advance \ycoord by -\ytop                              
\global\read\labeLfile to\justx                                       
\global\read\labeLfile to\justy                                       
\global\read\labeLfile to\labeL
\global\setbox\labox=\hbox{\labeL\hskip-0.3em}%    
\advance\vpos by-\ycoord                                              
\vskip-\vpos \vpos=\ycoord                                         
\hbox to\pswidth{\hskip\xcoord %                                 
\hbox to 0pt{\ifnum\justx>0\hss\fi%                                   
\vbox to0pt{%                                                         
\ifnum\justy<2\vss\fi%                                                
\copy\labox\kern0pt%  
\ifnum\justy>0\vss\fi}%                                               
\ifnum\justx<2\hss\fi}%                                               
\hss}%                                                                
\repeat%                                                              
\advance\vpos by-\psheight%                                           
\vskip-\vpos %                                                     
}\closein\labeLfile}
\def\figplace#1#2#3{
\openin\labeLfile=\figdir#1.lbl
\ifeof \labeLfile
       \immediate\write16{***Can't find \figdir#1.lbl; Skipping it.***}
\else  \closein\labeLfile
       \null\hskip#2\raise #3 \hbox{\figcrop{#1}}
\fi
}
\def\bbbone{{\mathchoice {\rm 1\mskip-4mu l} {\rm 1\mskip-4mu l}
{\rm 1\mskip-4.5mu l} {\rm 1\mskip-5mu l}}}

\def\NN{\mathbb{N}}

\newcommand{\bea} {\begin{eqnarray}}
\newcommand{\eea} {\end{eqnarray}}
\newcommand{\lp}  {\left(}
\newcommand{\rp}  {\right)}
\newcommand{\cF}{\mathcal F}
\newcommand{\cJ}{\mathcal J}
\newcommand{\cC}{\mathcal C}
\newcommand{\cN}{\mathcal N}
\newcommand{\cA}{\mathcal A}
\newcommand{\cZ}{\mathcal Z}
\newcommand{\cM}{\mathcal M}
\newcommand{\cT}{\mathcal T}

\newcommand{\ep}{\epsilon}
\newcommand{\si}{\sigma}
\newcommand{\ga}{\gamma}
\newcommand{\Ga}{\Gamma}
\newcommand{\al}{\alpha}

\newcommand{\De}{\Delta}
\newcommand{\ph}{\phi}

\def\Br{\overline}

\newtheorem{Theorem}{Theorem}[section]
\newtheorem{Lemma}[Theorem]{Lemma}
\newtheorem{Proposition}[Theorem]{Proposition}
\newtheorem{Corollary}[Theorem]{Corollary}
\newtheorem{Remark}[Theorem]{Remark}

\newtheorem{Definition}[Theorem]{Definition}
\newtheorem{Example}[Theorem]{Example}

\newcommand{\complex}{\mathbf C}
\newcommand{\Z}{\mathbf Z}
\newcommand{\bF}{\mathbb F}
\renewcommand{\P}{\mathbb P}        
\newcommand{\T}{\mathbb T}        

\renewcommand{\O}{\mathcal O} 
\newcommand{\I}{\mathcal I}

\newcommand{\R}{\mathcal R} 
\newcommand{\E}{\mathcal E} 
\newcommand{\F}{\mathcal F} 
\newcommand{\ux}{\mathbf x} 
\newcommand{\uy}{\mathbf y} 
 
\newcommand{\U}{\mathcal E} 
\newcommand{\Se}{\mathcal S} 
\newcommand{\Te}{\mathcal T} 
\newcommand{\ra}{\rightarrow}

\newcommand{\lra}{\longrightarrow}

\renewcommand{\ker}{\text{ker}\,}

\newcommand{\demo}{\noindent {\sc Proof.}\;}
 
\begin{document}
\title[Brill-Gordan loci and the Foulkes conjecture]{Brill--Gordan Loci, 
transvectants \\ and an analogue of the Foulkes conjecture} 
\author[Abdesselam and Chipalkatti]{Abdelmalek Abdesselam and Jaydeep Chipalkatti} 
\maketitle 

\parbox{12cm}{\small 
{\sc Abstract.} 
The hypersurfaces of degree $d$ in the projective space $\P^n$ correspond 
to points of $\P^N$, where $N = \binom{n+d}{d}-1$. Now assume $d=2e$ is 
even, and let $X_{(n,d)} \subseteq \P^N$ denote the subvariety of 
two $e$-fold hyperplanes.
We exhibit an upper bound on the Castelnuovo regularity of the 
ideal of $X_{(n,d)}$, and show that this variety is $r$-normal 
for $r \ge 2$. The latter result is representation-theoretic, and 
says that a certain $GL_{n+1}$-equivariant moprhism 
\[S_r(S_{2e}(\complex^{n+1}))\longrightarrow S_2(S_{re}(\complex^{n+1})) \]
is surjective for $r\ge 2$ ; a statement which is reminiscent 
of the Foulkes-Howe 
conjecture. For its proof, we reduce the statement to the case $n=1$, and 
then show that certain transvectants of binary forms are nonzero. 
The latter part uses explicit calculations with Feynman diagrams and 
hypergeometric series. 
For ternary quartics and binary $d$-ics, we give explicit generators for the 
defining ideal of $X_{(n,d)}$ expressed in the language of classical invariant theory.}

\vspace{5mm} 

\mbox{\small AMS subject classification (2000): 05A15, 14F17, 20G05, 81T18.}

\medskip 

\parbox{12cm} 

\medskip 

\section{Introduction} 

This article is addressed to a rather diverse audience: representation theorists,
algebraic geometers, combinatorialists, specialists in hypergeometric
series and angular momentum, as well as theoretical physicists working on 
quantum gravity using spin networks. In this introduction, we try to 
describe our results in a manner accessible to all.

\subsection{The Foulkes-Howe conjecture} 
One of the major problems in the representation theory of the general linear group 
is understanding the composition of Schur functors, variously known as 
plethysm or `external product' of symmetric functions. 
Even in the `simple' case of a composition of symmetric
powers $S_r(S_m(\complex^{n+1}))$ (which is the space of homogeneous polynomials
of degree $r$ in the coefficients of a generic homogeneous polynomial
of degree $m$ in $n+1$ variables), 
very little is known about its decomposition into irreducible representations
of $SL_{n+1}$. While trying to shed light on this problem, R.~Howe~\cite{Howe} 
constructed a natural equivariant map 
\[ S_r(S_m(\complex^{n+1}))\longrightarrow S_m(S_r(\complex^{n+1})). \] 
He conjectured that the map is injective if $r\le m$, and surjective if $r \ge m$, 
thereby giving a more precise form to a question raised by H.~O.~Foulkes~\cite{Foulkes}. 
(See~\cite{Briand,Brion2,Doran} for recent results and further references.)
More generally, for any integer $e \ge 1$, there is an
equivariant map 
\begin{equation}
 S_r(S_{me}(\complex^{n+1}))\longrightarrow S_m(S_{re}(\complex^{n+1})), 
\label{general.howemap} \end{equation}
which reduces to Howe's map for $e=1$. 
(An explicit definition of the map will be given in Section \ref{transvectants}.) 
An immediate question is whether this more general map also is surjective when $r \ge m$.
Our main result says that this is so for $m=2$. 

\begin{Theorem} \sl
The map 
\[ \alpha_r: S_r(S_{2e}(\complex^{n+1}))\longrightarrow 
S_2(S_{re}(\complex^{n+1})) \] 
is surjective for $r\ge 2$. 
\label{surj.alpha_r} \end{Theorem}
\begin{Remark} \rm 
The following result was recently proved by Rebecca Vessenes in her 
thesis (see~\cite[Theorem 1]{Vessenes}):
{\sl For any partition $\lambda$ and $r \ge 2$, the multiplicity of 
the irreducible Schur module $S_\lambda(\complex^{n+1})$ 
in $S_r(S_{2e}(\complex^{n+1}))$ 
is at least equal to its multiplicity in $S_2(S_{re}(\complex^{n+1}))$.} 
The theorem above of course implies this. The technique of tableaux counting 
used by her gives a similar (but slightly weaker) result 
(see [loc.~cit., Theorem 2]): 
{\sl For $r \ge 3$, any module $S_\lambda(\complex^{n+1})$ which has positive 
multiplicity in $S_r(S_{3e}(\complex^{n+1}))$ also has positive 
multiplicity in $S_3(S_{re}(\complex^{n+1}))$.} This is inaccessible by our 
method as it stands. 
\end{Remark} 
\begin{Remark} \rm To the best of our knowledge, the map~(\ref{general.howemap})
is first considered by Brion (see~\cite[\S 1.3]{Brion1}). 
He shows that there exists a constant 
$C(m,e,n)$, such that~(\ref{general.howemap}) is surjective for
$r \ge C(m,e,n)$. 
\end{Remark}
\subsection{Brill-Gordan loci}
In fact, we discovered Theorem \ref{surj.alpha_r} in the course 
of an entirely different line of inquiry. The context is as follows: 

The set of hypersurfaces of degree $d$ in $\P^n$ is parametrized by the 
projective space $\P^N$, where $N = \binom{n+d}{d}-1$. Assume that $d$ is 
even (say $d=2e$), and consider the subset of hypersurfaces 
which consist of two (possibly coincident) $e$-fold hyperplanes. 
In algebraic terms, we regard $\P^N$ as the space of degree $d$ forms in 
$n+1$ variables (up to scalars), and consider the set 
\[ \{ [F] \in \P^N: F = (L_1\, L_2)^e \; \; \text{for some 
linear forms $L_1,L_2$} \}. \] 
This is a projective subvariety of $\P^N$, which we denote by $X_{(n,d)}$. 
Throughout we exclude the trivial case $n=1,d=2$, and write 
$X$ for $X_{(n,d)}$ if no confusion is likely. 
The imbedding $X \subseteq \P^N$ is stable for the natural action 
on $SL_{n+1}$. 

This construction is modelled after the variety of {\sl 
totally decomposable forms}, defined as 
\[ Y = \{ [F] \in \P^N: F = L_1 L_2 \dots L_d \; \; \text{for some 
linear forms $L_i$} \}. \] 
Brill \cite{Brill1,Brill2} and Gordan \cite{Gordan} considered the problem of 
finding $SL_{n+1}$-invariant defining equations for $Y$. In classical terms, 
we are to find a set of concomitants of a generic $(n+1)$-ary $d$-ic $F$ 
which vanish iff $F$ belongs to $Y$. (It turns out that there exists 
such a set of concomitants in degree $d+1$; see~\cite[Ch.~4]{GKZ} for a 
modern account of Brill's work.) 
Due to this obvious analogy, we may call $X$ a Brill-Gordan locus. 

This project began as an attempt to find defining equations for $X$. 
This led to the following statement about the homogeneous ideal $I_X$. 

\begin{Theorem}[Main Theorem] \sl 
The ideal $I_X$ is $m_0$-regular with $m_0 = \lceil 2n+1 - \frac{n}{e} \rceil$. 
{\it A fortiori}, $X$ is scheme-theoretically defined by equations of 
degree at most $m_0$. 
\label{main.theorem} 
\end{Theorem}
In order to prove the first statement, it is necessary to show that the cohomology groups 
$H^i(\P^N,\I_X(m_0-i))$ are zero for $i \ge 1$. The case $i=1$ is the hardest part of 
the proof. It follows once we show that the morphism 
\[ H^0(\P^N, \O_{\P^N}(m_0-1)) \lra H^0(\P^N, \O_X(m_0-1)) \] 
is surjective. Once both sides are identified \emph{qua} 
$SL_{n+1}$-representations, we are reduced to showing that the morphism 
\[ \alpha_r: S_r(S_{2e}(\complex^{n+1}))\longrightarrow S_2(S_{re}(\complex^{n+1})) \] 
is surjective for $r = m_0-1$. While attempting to prove this, we realized that 
the surjectivity in fact holds for {\sl all} $r \ge 2$, which is 
Theorem~\ref{surj.alpha_r}. Alternately said, 
the variety $X$ is $r$-normal for $r \ge 2$. Since $\alpha_2$ is an isomorphism, 
$I_X$ contains no degree $2$ forms. 
\subsection{Examples}
Notice that $m_0=3$ when $n=1$, hence $X$ is defined by cubic equations in this 
case. We describe these equations explicitly in section \ref{binarycase}.
The answer is formulated in terms of degree $3$ covariants of binary $d$-ics 
(in the sense of \cite{GrYo}). To wit, we exhibit a finite set of covariants 
$\{\Phi_i\}$ such that 
{\sl a binary $d$-ic $F$ belongs to $X_{(1,d)}$, iff $\Phi_i(F)=0$ for all $i$.} 

The example of ternary quartics (i.e., the case $n=2,d=4$) is worked out in 
section \ref{ternaryquartics}. It partly relies upon some elimination-theoretic 
computations done in Macaulay-2. 

\begin{Remark} \rm 
One can define a Brill-Gordan locus rather generally. 
Associated to a partition $\lambda = (\lambda_1,\lambda_2,\dots)$ of $d$, 
we have a subvariety $X^\lambda \subseteq \P^N$ of 
forms which factor as $F = \prod L_i^{\lambda_i}$. It is 
a natural problem to find $SL_{n+1}$-invariant equations for 
this variety. The case $\lambda=(d)$ corresponds to the Veronese imbedding
(see \cite{JoeH}), and $\lambda=(1^d)$ is the case 
considered by Brill and Gordan. Alternately, $X^{(1^d)}$ can be 
identified with the variety of Chow forms of degree $d$ zero-cycles in $\P^n$ 
(see \cite{GKZ}). 
A result for the case $\lambda =(\lambda_1,\lambda_2)$ with 
$\lambda_1 > \lambda_2$ is in preparation. 

In the case of binary forms, $X^\lambda$ is the so called `coincident root locus'. 
It was first studied by Cayley in~\cite{Cayley1}, and has received 
recent attention in the work of Dixmier~\cite{Dixmier1,Dixmier2},
Weyman~\cite{Weyman1,Weyman2,Weyman3} and the second 
author~\cite{Chipalkatti1,Chipalkatti2}.
There are also connections to singularity theory~\cite{Feher}, and the combinatorics of
integrable systems~\cite{Kasatani}. A set of $SL_2$-invariant defining 
equations is known for binary forms (see ~\cite{Chipalkatti2}), however 
the ideal $I_X$ is not well-understood.

Of course we can reformulate the question by allowing factors of higher degree. 
For instance, the quartic plane curves which split into a line and a cubic fill up 
an $11$-dimensional subvariety of $\P^{14}$; we do not know its defining equations. 
\end{Remark} 

\subsection{The Proof of Theorem \ref{surj.alpha_r}} 
We give a short description of the principal steps in the proof. 
By a formal argument, it suffices to consider the case $n=1$. 
Now we have a plethysm decomposition 
\begin{equation} S_2(S_{re}(\complex^2)) = 
\bigoplus\limits_p S_{rd -4p}(\complex^2), 
\end{equation}
where the direct sum is quantified over $0 \le p \le 
\lfloor \frac{re}{2} \rfloor$. Let $\pi_p$ denote the 
projection onto the $p$-th summand. By Schur's lemma, it is enough to 
show the following: 
\begin{Proposition}[Key Proposition] \sl 
When $n=1$, $r\ge 2$ and $0 \le p \le 
\lfloor \frac{re}{2} \rfloor$,
the morphism $\pi_p \circ \alpha_r$ is nonzero. 
\label{key.prop} \end{Proposition}

The proof is by induction on $r$, and occupies the bulk of the paper. 
The initial result for $r=2$, and the induction step are respectively 
proved in Lemmata \ref{lemma.A} and \ref{lemma.B}. 
In either case, the crux of the result consists in showing that 
certain Clebsch-Gordan coefficients (or what are the essentially the same, 
Wigner's $3j$-symbols) are nonzero. 

We give two proofs of Lemma \ref{lemma.A}.
The first is a combinatorially explicit calculation with Feynman diagrams
(used here as the pictorial counterpart of classical covariants)
which explains {\em why} the corresponding coefficient is nonzero.
The second is perhaps less transparent, but it
allows a closed form evaluation, thanks to Dixon's summation theorem
for the ${}_3 \bF_2$ hypergeometric series.
The proof of Lemma~\ref{lemma.B} uses Feynman diagrammatic generating
function techniques. These are implicit in the work of 
J.~Schwinger~\cite{Schwinger} (which is based on the second quantization formalism), 
and its restatement by V.~Bargmann~\cite{Bargmann} which uses Gaussian integration.
At this point, the use of analysis (akin to Bargmann's) would be a tempting shortcut. 
However, this would have obscured the fact that what is at play is purely 
combinatorial algebra; and, except as a guiding principle, there is no real need 
for transcendental methods. 

\begin{Remark} \rm 
This is an instance of the combinatorial underpinnings behind the 
invariant theory of binary forms. The latter is a fascinating subject 
(see~\cite{Elliott,Glenn,GrYo,Salmon} for classical accounts), with 
ramifications in many fields of current mathematical and physical interest. 
For instance, it makes an appearance in the quantum theory of angular 
momentum~\cite{Biedenharn1, Biedenharn2}, classical hypergeometric 
series~\cite{Gustafson}, the spin network approach to quantum 
gravity~\cite{Penrose, Rovelli}, as well as knot and $3$-manifold 
invariants~\cite{Carter}. Modern presentations of the classical invariant theory 
of binary forms may be found in \cite{KungRota} and \cite{Olver}. 
\end{Remark} 

\subsection{Transvectants} 
We will rephrase Proposition~\ref{key.prop} as a statement about transvectants of 
binary forms. We begin by recalling the latter notion. 

Let $A(x_0,x_1)$ and $B(x_0,x_1)$ be binary forms of degrees $a,b$. 
Introduce new variables $(y_0,y_1)$, and consider the differential operator 
\[ \Omega = 
\frac{\partial^2}{\partial x_0 \, \partial y_1}
- \frac{\partial^2}{\partial x_1 \, \partial y_0}, 
\] 
usually known as Cayley's Omega operator. If $k$ is a nonnegative 
integer, then the $k$-th transvectant of $A,B$ is defined to be 
\begin{equation} (A,B)_k = 
\frac{(a-k)! \, (b-k)!}{a! \, b!} \, [ \, \Omega^k \, 
A(x_0,x_1) \, B(y_0,y_1)]_{\uy:=\ux}. 
\label{trans.defn} \end{equation}
(This is interpreted as follows: change $(x_0,x_1)$ to $(y_0,y_1)$ in $B$, apply 
$\Omega$ in all $k$-times to the product $A \, B$, and finally substitute $x_i$ for $y_i$.) 
By construction, $(A,B)_k$ is a binary form of degree $a+b-2k$. It is identically 
zero if $k > \min \{a,b\}$. A general account of 
transvectants may be found in \cite{Glenn,GrYo,Olver}. 

Now the Key Proposition is equivalent to the following statement:
\begin{Proposition} \sl 
Let $Q$ be a generic binary form of degree $r\ge 2$. 
Then, for any integers $e,p$ such that $e \ge 1,0\le\ p\le \lfloor \frac{re}{2} \rfloor$,
the transvectant $(Q^e,Q^e)_{2p}$ is not identically zero. 
\label{key.transvect}
\end{Proposition}
A proof of the equivalence is given in section \ref{prop.eq}. 
\begin{Example} \rm 
In general, it may be a nontrivial matter to show that a given 
transvectant expression is (or is not) identically zero. As a simple exercise, the 
reader should check that $(F,(F,F)_2)_5 = 0$ for any binary quintic $F$.  
We will later see that the {\sl odd} transvectants $(Q^e,Q^e)_{2p+1}$ are zero. 
\end{Example} 

\begin{Remark} \rm 
With some more work (which we do not do), one can 
probably trace through our inductive proof of
the Key Proposition in order to obtain an explicit formula for the
transvectant $(Q^e,Q^e)_{2p}$. The latter is an $SL_2$-invariant 
function of $r+1$ points on the Riemann sphere $\P^1$: namely the roots 
of $Q$, and the point with homogeneous coordinates $x_0,x_1$. At 
least on a very intuitive level, our induction on $r$ can be thought of as 
degenerating the surface in order to separate the points to the extent allowed 
by the stability criterion for ${\Br\cM}_{0, r+1}$ 
(the moduli space of genus zero curves with $r+1$ labelled points). 

It would be intriguing if one could make this intuition precise
using the powerful method of equivariant localization on the 
corresponding strata of ${\Br\cM}_{0, r+1}$. This would open the
door to the application of these techniques to the calculation of 
new formulae for other specific covariants of binary forms.
\end{Remark} 

\subsection{Symmetric functions} 
If we express the previous transvectant in terms of the roots
of $Q$ and then dehomogenize, this becomes a nonvanishing statement
for ordinary symmetric functions defined as sums over magic squares or
transportation matrices with integer entries. (See~\cite{Diaconis}
for a recent review on the fascinating combinatorics of these objects.)

We start with $r+1$ variables, $z_1,\ldots,z_r$ and $t$. 
Let $\cM$ denote the set of $(r+1) \times (r+1)$ 
matrices $M=(m_{ij})_{1\le i,j\le r+1}$ satisfying the following 
conditions: 
\begin{itemize} 
\item 
all $m_{i,j}$ are nonnegative and the diagonal entries are zero, 
\item 
the row and column sums are given by the vector \newline 
$(\underbrace{e,\ldots,e}_r, re-2p)$. 
\end{itemize} 
Now define the following symmetric function in the $z_i$, with $t$ 
as a parameter. 
\begin{equation}
\begin{aligned} {} & \cT_{r,e,p}(t; z_1,\ldots,z_r) =  \\
& \sum_{M\in \cM}
\frac{\prod\limits_{1\le i,j\le r} (z_i-z_j)^{m_{ij}}
\prod\limits_{1\le i\le r} (t-z_i)^{m_{i,r+1}}
\prod\limits_{1\le j\le r} (t-z_j)^{m_{r+1,j}}}
{\prod\limits_{1\le i,j\le r+1}m_{ij}!} 
\end{aligned} \label{tau} \end{equation}
We now have the following result. 
\begin{Proposition} \sl 
For any $r\ge 2$, $e\ge 1$ and $0\le\ p\le \lfloor \frac{re}{2} \rfloor$,
the function 
$\cT_{r,e,p}(t; z_1,\ldots,z_r)$ does not identically vanish. 
\label{magic.squares}
\end{Proposition} 
The case when $re$ is even and $p=\frac{re}{2}$
does not involve $t$ and is perhaps the most aesthetically pleasing: 
it reduces to a sum over $r \times r$ magic squares, with row and 
column sums given by $e$. It would be an interesting problem to express 
$\cT_{r,e,p}$ in terms of Schur functions.
\subsection{The symbolic method} \label{symbolic.method}
We will freely use the symbolic method of classical invariant theory 
(see~\cite{GrYo, Olver}). Since this has ceased to be a part of the algebraists'
standard repertoire, a few words of explanation are in order. 
The symbolic method is a very powerful tool, with a simple underlying principle. 

As an example, take four pairs of binary variables
$a=(a_0,a_1)$, $b=(b_0,b_1)$, $c=(c_0,c_1)$ and $d=(d_0,d_1)$. 
Let the {\sl symbolic bracket}  $(a \, b)$ stand for $a_0 \, b_1- b_0 \, a_1$ etc. 
Now consider the following algebraic expression 
\begin{equation}
E = (a\, b)^2 \, (c \, d)^2 \, (a\, c) \, (b\, d). 
\label{discriminant}
\end{equation}
Each letter occurs three times, hence classically $E$ represents an 
invariant of binary cubics. This is interpreted as follows: 
if $F(x_0,x_1)$ denotes the generic binary cubic, then $E$ represents the 
algebraic expression obtained by applying the differential operator 
\[
(\frac{1}{3!})^4 \, F(\frac{\partial}{\partial a_0},\frac{\partial}{\partial a_1}) \, 
F(\frac{\partial}{\partial b_0},\frac{\partial}{\partial b_1}) \, 
F(\frac{\partial}{\partial c_0},\frac{\partial}{\partial c_1}) \, 
F(\frac{\partial}{\partial d_0},\frac{\partial}{\partial d_1}) \] 
to $(a\, b)^2(c \, d)^2(a\, c)(b\, d)$. 
The result is a homogeneous degree $4$ polynomial 
in the coefficients of $F$. (Up to a scalar, it is the discriminant of $F$.)
This interpretation is the reverse or dualized form of the one given 
in~\cite[Appendix I]{GrYo}.
We believe that it offers several advantages in simplicity 
and flexibility: for instance the possibility 
of iteration, or mixed interpretation (where some variables 
are taken as `actual' and others as symbolic within the same computation). 
Symbolic letters are nothing more than auxiliary variables which 
are differentiated out in the final interpretation of the expressions
at hand. 

The symbolic method has a rather close resemblance to modern calculational methods 
from physics~(e.g., see \cite{Cvitanovic}). 
The formal brackets can be seen as the result of differentiating (or integrating)
out anticommuting Fermionic variables (see~\cite{Clifford}). 
In such calculations one often quickly faces an inflation of the number of 
letters needed, and may wonder how to label them. Perhaps one can do this with 
points of an infinite variety like a string worldsheet, thereby
organizing the collection of these variables into a `quantum field'.
This suggests the question of interpreting topological field theoretic constructions 
along these lines--compare~\cite{Cattaneo} and \cite{Kontsevich}.
%%%%%%%%%%%%%%%%%%%%%%%%%%%%%%%%%%%%%%%%%%%%%%%%%%%%%%%%%%%%%%%%%%%
\section{Preliminaries} 
In this section we establish the set-up and notation. 
All terminology from algebraic geometry follows \cite{Ha}. 

The base field will be $\complex$. 
Let $V$ denote a complex vector space of dimension $n+1$, and 
write $W=V^*$. If $\lambda$ is a partition, then $S_\lambda(-)$ will 
denote the associated Schur functor. In particular, $S_d(-)$ denotes the 
symmetric power. All subsequent constructions will be $SL(V)$-equivariant;
see \cite[Ch.~6 and 15]{FH} for the relevant representation theory. 
Normally we suppress the reference to $V$ whenever it is understood 
from context. Thus, for instance, 
$S_r(S_d)$ stands for $\text{Sym}^r(\text{Sym}^d \, V)$.

Fix a positive integer $d=2e$, and let $N = \binom{n+d}{d}-1$. 
Given the symmetric algebra 
\[ R = \bigoplus_{r \ge 0} S_r(S_d \, V), \] 
the space of degree $d$ hypersurfaces in $\P V$ is identified 
with 
\[ \P^N = \P \, S_d \, W = \text{Proj} \; R. \] 
Now define 
\begin{equation}
X_{(n,d)} = \{ [F] \in \P^N: F = (L_1 L_2)^e \;\; \text{for 
some $L_1,L_2 \in W$} \}. 
\end{equation}
This is an irreducible $2n$-dimensional projective subvariety of 
$\P^N$. 

\smallskip 

Recall the definition of regularity according to Mumford 
\cite[Ch.~6]{Mum}. 
\begin{Definition} \sl 
Let $\F$ be a coherent $\O_{\P^N}$-module, and $m$ an integer. 
Then $\F$ is said to be 
$m$-regular if $H^q(\P^N,\F(m-q)) = 0$ for $q \ge 1$. 
\end{Definition} 
It is known that $m$-regularity implies $m'$-regularity for 
all $m' \ge m$. 
Let $M$ be a graded $R$-module containing no submodules of 
finite length. Then (for the present purpose) we will say that $M$ is 
$m$-regular if its sheafification ${\widetilde M}$ is. In our case, 
$M=I_X$ (the saturated ideal of $X$), and ${\widetilde I_X} = \I_X$. 

We have the usual short exact sequence 
\begin{equation} 
0 \ra \I_X \ra \O_{\P^N} \ra \O_X \ra 0. 
\label{ses1} \end{equation} 
The map 
\begin{equation} \P W \times \P W \stackrel{f}{\lra} \P S_d \, W, \quad 
(L_1,L_2) \lra (L_1L_2)^e
\end{equation}
induces a natural isomorphism of $X$ with the quotient 
$(\P W \times \P W)// \Z_2$, and of the structure sheaf $\O_X$ with 
$(f_* \O_{\P W \times \P W})^{\Z_2}$. 

Using the Leray spectral sequence and the K{\"u}nneth formula, 
\[ \begin{aligned} 
{} & H^q(\P^N, f_* \O_{\P W \times \P W}(r)) = \\ 
\bigoplus\limits_{i+j=q} & 
H^i(\P W, \O_{\P^n}(re)) \otimes 
H^j(\P W, \O_{\P^n}(re)). 
\end{aligned} \] 
This group can be nonzero only in two cases: $i,j$ are either both 
$0$ or both $n$ (see \cite[Ch.~III,\S 5]{Ha}). 
Now $H^q(\O_X(r))$ is the $\Z_2$ invariant part of 
$H^q(f_* \O_{\P W \times \P W}(r))$ for any $q$, which gives the 
following corollary. 

\begin{Corollary} \sl 
We have an isomorphism $H^0(\O_X(r)) = S_2(S_{re})$ for 
$r \ge 0$. Moreover $H^{2n}(\O_X(r)) =0$ for $re \ge -n$. 
If $q \neq 0$ or $2n$, then $H^q(\O_X(r))=0$. 
\qed \end{Corollary} 
%%%%%%%%%%%%%%%%%%%%%%%%%%%%%%%%%%%%%%%%%%%%%%%%%%
\section{The Proof of the Main Theorem} 
In this section we begin the proof of Theorem \ref{main.theorem}. 
Modulo some cohomological arguments, it will reduce to the statement of Theorem 
\ref{surj.alpha_r}. The latter will be proved in sections \ref{transvectants} and 
\ref{Feynman}. 

Define the predicate 
\[ \R(q): \; H^q(\P^N, \I_X(m_0-q)) = 0. \] 
We would like to show $\R(q)$ for $q \ge 1$. 
Tensor the short exact sequence (\ref{ses1}) by $\O_{\P^N}(m_0-q)$, and  
consider the piece 
\begin{equation}  \dots \ra H^{q-1}(\O_X(m_0-q)) \ra H^q(\I_X(m_0-q)) \ra 
          H^q(\O_{\P^N}(m_0-q)) \ra \dots 
\end{equation}
from the long exact sequence in cohomology. 
We claim that if $q > 1$ then the first and third terms vanish, 
hence $\R(q)$ is true. This is clear if $q \neq 2n+1$. 
By the choice of $m_0$, we have 
\[ e \, (m_0-2n-1) \ge -n, \] 
implying that $H^{2n}(\O_X(m_0-2n-1))=0$. Hence the claim is still true if 
$q = 2n+1$. 
It remains to prove $\R(1)$, which is the special case $r=m_0-1$ of the following result. 
\begin{Proposition} \sl 
Let $r \ge 2$. Then the morphism 
\[ \alpha_r: H^0(\O_{\P^N}(r)) \lra H^0(\O_X(r)) \] 
is surjective. 
\end{Proposition}
\demo 
The map $f$ can be factored as 
\[ \P W \times \P W \lra 
\P S_e \, W \times \P S_e \, W \lra \P S_d \, W. \] 
Tracing this backwards, we see that $\alpha_r$ is the composite 
\begin{equation}
S_r(S_d) \stackrel{1}{\lra} S_r(S_e \otimes S_e) 
\stackrel{2}{\lra} S_r(S_e) \otimes S_r(S_e) 
\stackrel{3}{\lra} S_{re} \otimes S_{re} 
\stackrel{4}{\lra} S_2(S_{re}), 
\label{alpha.r} \end{equation}
where 1 is given by applying $S_r(-)$ to the coproduct map, 2 is the  
projection coming from the `Cauchy decomposition' (see \cite{Akin}), 3 is the 
multiplication map, and 4 is the symmetrisation. 
Now we have a plethysm decomposition 
\begin{equation} H^0(\O_X(r)) = S_2(S_{re}) = 
\bigoplus\limits_p S_{(rd-2p,2p)}, 
\label{oxr} \end{equation}
where the direct sum is quantified over $0 \le p \le 
\lfloor \frac{re}{2} \rfloor$. Let $\pi_p$ denote the 
projection onto the $p$-th summand. Since any finite dimensional 
$SL(V)$-module is completely reducible, the cokernel of $\alpha_r$ is a 
direct summand of $H^0(\O_X(r))$. We will show that 
$\pi_p \circ \alpha_r \neq 0$ for any $p$, then Schur's lemma will 
imply that the cokernel is zero. 

The entire construction is functorial in $V$, hence if 
$U \subseteq V$ is any subspace, then the diagram 
\[ \diagram
S_r(S_d \, U) \dto \rto & S_{(rd-2p,2p)} \, U \dto \\
S_r(S_d \, V) \rto & S_{(rd-2p,2p)} \, V 
\enddiagram \] 
is commutative. The vertical map on the left is injective. 
If we further assume that $\dim U \ge 2$, then the vertical map on the right is 
injective as well. (Recall that $S_\lambda(V)$ vanishes if and only if the 
number of parts in $\lambda$ exceeds $\dim V$.) Hence we may as well assume that 
$\dim V=2$. Thus we are reduced to the Key Proposition~\ref{key.prop} (see the 
Introduction); we defer its proof to Sections~\ref{transvectants} and~\ref{Feynman}. 
\qed 

This reduction argument can be understood as follows: $\pi_p \circ \alpha_r$ is 
a formal multilinear construction involving $n+1$ variables. If it gives 
a nonzero result when all but two of the variables are set to zero, then it must 
have been nonzero to begin with.

Note the following simple corollary to the Main Theorem. 
\begin{Corollary} \sl 
In the Grothendieck ring of finite-dimensional $SL(V)$-modules, we 
have the equality
\[ [(I_X)_r] = 
[S_r(S_d)]- \sum\limits_{0 \le p \le \lfloor \frac{re}{2} \rfloor} 
[S_{(rd-2p,2p)}] \] 
Here $[-]$ denotes the formal character of a representation. 
\label{gr} \end{Corollary} 
\demo This follows because $(I_X)_r = \ker \alpha_r$.  \qed 

\smallskip 

Decomposing the plethysm $S_r(S_d)$ into irreducible submodules is 
in general a difficult problem. Explicit formulae are known only in very 
special cases -- see \cite{Chen,Macdonald} and the references therein. 
In particular the decomposition of $S_3(S_d)$ is given by Thrall's 
formula (see \cite{Plunkett1}), and then $(I_X)_3$ can be 
calculated in any specific case. 
Note that $(I_X)_2 =0$, i.e., there are no quadratic polynomials 
vanishing on $X$. 
%%%%%%%%%%%%%%%%%%%%%%%%%%%%%%%%%%%%%%%%%%%%
\section{Ternary quartics} \label{ternaryquartics} 
Assume $n=2,d=4$. We will identify the 
generators of $I_X$ as concomitants of ternary quartics in the sense 
of classical invariant theory. We will partly rely upon some computations done 
using the program Macaulay-2. 

By the Main Theorem we know that the generators of 
$I_X$ lie in degrees $\le 4$. We will find them using an 
elimination theoretic computation. Define 
\[ \begin{aligned} 
   L_1 & = a_0 \, x_0 + a_1 \, x_1 + a_2 \, x_2, \quad 
   L_2 = b_0 \, x_0 + b_1 \, x_1 + b_2 \, x_2, \\
   F & = c_0 \, x_0^4 + c_1 \, x_0^3 \, x_1 + \dots + c_{14} \, x_2^4;
\end{aligned} \] 
where the $a,b,c$ are indeterminates. 
Write $F = (L_1L_2)^2$ and equate the coefficients of the monomials 
in $x_0,x_1,x_2$. This expresses 
each $c_i$ as a function of $a_0,\dots,b_2$, and hence defines a ring map 
\[ \complex \, [c_0,\dots,c_{14}] \lra \complex \, [a_0, \dots,b_2]. \] 
The kernel of this map is $I_X$. When we calculated it using Macaulay-2, it 
turned out that in fact all the minimal generators are in degree $3$, 
hence it is enough to look at the piece $(I_X)_3$. 
By Corollary \ref{gr} and Thrall's formula, 
\[ (I_X)_3 = S_{(9,3)} \oplus S_{(6,0)} \oplus S_{(6,3)} 
\oplus S_{(4,2)} \oplus S_{(0,0)}. \] 
Note that an inclusion 
\[ S_{(m+n,n)} \subseteq (I_X)_3 \subseteq S_3(S_4) \] 
corresponds to a concomitant of ternary quartics having degree $3$,
order $m$ and class $n$. (This correspondence is fully explained in 
\cite{Chipalkatti3}.) 
For instance, $S_{(9,3)}$ corresponds to a concomitant of 
degree $3$, order $6$ and class $3$. 

It is not difficult to identify the concomitants symbolically 
(see [loc.~cit.] for the procedure). In our case, the 
summands respectively correspond to 
\begin{equation} \begin{array}{ll}
\alpha_x^2 \, \beta_x^3 \, \gamma_x \, (\alpha \, \gamma \, u)^2 
(\beta \, \gamma \, u), & 
\alpha_x^2 \, \beta_x^2 \, \gamma_x^2 \, (\alpha \, \beta \, \gamma)^2, \\ 
\alpha_x^2 \, \beta_x \, (\beta \, \gamma \, u)^2 
(\alpha \, \gamma \, u) (\alpha \, \beta \, \gamma), & 
\alpha_x \, \beta_x \, (\alpha \, \gamma \, u) (\beta \, \gamma \, u) 
(\alpha \, \beta \, \gamma)^2, \\ 
(\alpha \, \beta \, \gamma)^4. 
\end{array} \label{x2} \end{equation}
We can rephrase the outcome in geometric terms: 
\begin{Theorem}\sl 
Let $F$ be a ternary quartic with zero scheme $C \subseteq \P^2$. 
Then $C$ consists of two (possibly coincident) double lines 
iff all the concomitants in (\ref{x2}) vanish on $F$. 
\end{Theorem}

A similar result for any $n=1$ and any (even) $d$ will be deduced in 
Section \ref{binarycase}. 
%%%%%%%%%%%%%%%%%%%%%%%%%%%%%%%%%%%%%%%%%%%%
\section{The proof of Proposition \ref{key.prop}} \label{transvectants} 
In this section we will break down Proposition~\ref{key.prop} into two 
separate questions about transvectants of binary forms. 

\subsection{} 
We begin by describing the map $\alpha_r$ from (\ref{alpha.r}) in coordinates. 
(It is as yet unnecessary to assume $\dim V =2$.) 
Let 
\[ \ux^{(i)} = (x_0^{(i)}, \dots, x_n^{(i)}), \quad 
1 \le i \le r, \] 
be $r$ sets of $n+1$ variables. We will also introduce one set of their `copies' 
\[ \uy^{(i)} = (y_0^{(i)}, \dots, y_n^{(i)}), \quad 
1 \le i \le r. \] 

Let $F_i(\ux^{(i)}), 1 \le i \le r$ be degree $d$ forms, then the image 
$\alpha_r(\bigotimes\limits_{i=1}^r F_i)$ is calculated as follows: 
\begin{itemize} 
\item 
Apply the polarization operator 
\[ \sum\limits_{\ell=0}^n y^{(i)}_\ell 
\frac{\partial}{\partial x^{(i)}_\ell} \] 
to each $F_i$ altogether $e$ times, and denote the result by 
$F_i(\ux^{(i)},\uy^{(i)})$. 
\item 
Take the product $\prod\limits_i F_i(\ux^{(i)},\uy^{(i)})$, and make 
substitutions 
\[ x^{(i)}_\ell = x_\ell, \quad y^{(i)}_\ell = y_\ell, \] 
for all $i,\ell$. (This is tantamount to `erasing' the upper indices.)
This gives a form having degree $re$ in each set $\ux,\uy$. Since it is symmetric in 
$\ux,\uy$, we may think of it as an element of $S_2(S_{re})$. It is  
the image of $\otimes F_i$ via $\alpha_r$. 
\end{itemize} 

\begin{Remark} \rm 
The map 
\[  S_r(S_{me}(\complex^{n+1}))\longrightarrow S_m(S_{re}(\complex^{n+1})), \] 
mentioned in the introduction is constructed similarly. That is, 
we introduce $m-1$ sets of copies $\uy^{(i)}, \dots, {\mathbf q}^{(i)}$, 
and then apply $e$ times the polarization operator 
\[ \left(\sum\limits_{\ell=0}^n y^{(i)}_\ell \frac{\partial}{\partial x^{(i)}_\ell}
\right) 
\dots \left(\sum\limits_{\ell=0}^n q^{(i)}_\ell \frac{\partial}{\partial x^{(i)}_\ell}
\right) \] 
to each degree $me$ form $F_i$. 
\end{Remark} 

\subsection{} \label{prop.eq}
Suppose now that $\dim V = 2$. Given a $G(\ux,\uy) \in S_2(S_{re})$, the 
form $\pi_p(G)$ is obtained, up to a nonzero numerical multiple,
by calculating $\Omega^{2p} \, G$, and 
setting $\uy=\ux$. We will now show that 
Proposition \ref{key.prop} is equivalent to Proposition~\ref{key.transvect}. 

\smallskip 

\demo Let us write symbolically $F_i(\ux) = ({h^{(i)}_\ux})^d$, where 
\[ h^{(i)}_\ux = h_{i,0} \, x_0 + h_{i,1} \, x_1 \] 
are linear forms. Then, following the recipe of the previous section, 

\begin{equation} 
\pi_p \circ \alpha_r \, (\otimes F_i) = \underbrace{(Q^e,Q^e)_{2p}}_{\mathcal A}, 
\end{equation}
where $Q = \prod\limits_{i=1}^r h^{(i)}_\ux$.
The right hand side is to be interpreted as follows: we formally 
calculate ${\mathcal A}$ as a transvectant, and then substitute the
actual coefficients 
of $F_i$ for the monomials $h_{i,0}^{d-j} \, h_{i,1}^j$. 
By the discussion of Section~\ref{symbolic.method}, 
this amounts to applying the differential operator
\[
\frac{1}{(d!)^r} \, \prod_{i=1}^r \, F_i 
(\frac{\partial}{\partial h_{i,0}},
\frac{\partial}{\partial h_{i,1}} ) \]
to the polynomial ${\mathcal A}(\{h_{i,0},h_{i,1}\}_i,x_0,x_1)$. 

Now assume Proposition \ref{key.prop}. 
This implies that ${\mathcal A}$ as an algebraic function of the $\{h_{i,0},h_{i,1}\}$
is not identically zero. Hence it is possible to {\sl specialize} the $h$ to
some complex numbers so that ${\mathcal A}$ remains nonzero. This specializes $Q$ to a
binary $r$-ic for which $(Q^e,Q^e)_{2p} \neq 0$, which shows Proposition~\ref{key.transvect}. 

For the converse, assume the existence of a $Q$ such that the
transvectant above is nonzero. It factors 
as (say) $Q = \prod\limits_{i=1}^r l_i$. Then letting $F_i=l_i^d$, we see that 
$\pi_p \circ \alpha_r(\otimes F_i) \neq 0$. \qed 

\smallskip 

Now we will prove Proposition \ref{key.prop} by an induction on $r$. 
We reformulate $r=2$ as a separate lemma. 
\begin{Lemma} \sl 
If $Q$ is the generic binary quadratic, then 
$(Q^e,Q^e)_{2p} \neq 0$ for $0 \le p \le e$. \label{lemma.A} \end{Lemma}
\demo This will directly follow from formula (\ref{Lemmadixon1}) in Section \ref{Feynman}. 
\qed 

\smallskip 

\subsection{The induction step} 
For the transition from $r$ to $r+1$, consider the commutative diagram 
 
\[ \diagram 
S_r(S_d) \otimes S_d \dto \rto^{\alpha_r \otimes 1} & 
S_2(S_{re}) \otimes S_d \dto^{u_r} \\ 
S_{r+1}(S_d) \rto_{\alpha_{r+1}} & S_2(S_{re+e}) \enddiagram \] 

Assume that $\alpha_r$ (and hence $\alpha_r \otimes 1$) is surjective.
If we show that $u_r$ is 
surjective, then it will follow that $\alpha_{r+1}$ is 
surjective. We need to understand the action of $u_r$ on the 
summands of the decomposition (\ref{oxr}). The map 
\[ u_r^{(p,p')}: S_{rd-4p} \otimes S_d \lra S_{(r+1)d-4p'} \] 
is defined as the composite 
\[ \begin{aligned} 
{} & S_{rd-4p} \otimes S_d \ra (S_{re} \otimes S_{re}) \otimes S_d 
\ra (S_{re} \otimes S_{re}) \otimes (S_e \otimes S_e) \ra \\
& S_{(r+1)e} \otimes S_{(r+1)e} \ra S_{(r+1)d-4p'}. 
\end{aligned} \] 
Let $A \in S_{rd-4p}, B \in S_d$. We will follow the sequence of 
component maps and get a recipe for calculating the image 
$u_r^{(p,p')}(A \otimes B)$. Let 
\[ \Gamma_g = \sum \limits_{i=0}^{re} \binom{re}{i} \, 
g_i \, x_0^{re-i} x_1^i, \quad
\Gamma_h = \sum \limits_{i=0}^{re} \binom{re}{i} \, 
h_i \, x_0^{re-i} x_1^i, \]
be two {\sl generic} forms of degree $re$. That is to say, the 
$g_i,h_i$ are thought of as independent indeterminates. (Of course these 
$h_i$ are unrelated to the ones in the last section.) 
\begin{itemize} 
\item 
Let $T_1 = (\Gamma_g,\Gamma_h)_{2p}$ and 
$T_2 = (A,T_1)_{rd-4p}$. Then $T_2$ does not involve $x_0,x_1$. 
\item 
Obtain $T_3$ by making the substitutions 
\[ g_i = x_1^{re-i} (-x_0)^i, \quad 
   h_i = y_1^{re-i} (-y_0)^i, \] 
in $T_2$. 
\item 
Let 
\[ T_4 = (y_0 \frac{\partial}{\partial x_0} + 
          y_1 \frac{\partial}{\partial x_1})^{e} \, B, \] 
and $T_5 = T_3 \, T_4$. 
\item 
Let $T_6 = \Omega^{2p'} \, T_5$. 
Finally $u_r^{(p,p')}(A \otimes B)$ is obtained by substituting 
$x_0,x_1$ for $y_0,y_1$ in $T_6$. 
\end{itemize} 
Hence it is enough to show the following: 

For $p'$ in the range $0 \le p' \le \frac{(r+1)e}{2}$, there 
exists a $p$ such 
that $u_r^{(p,p')}(A \otimes B)$ is nonzero for some forms $A,B$ of 
degrees $rd-4p,d$ respectively. This will prove the surjectivity of 
$u_r$ and complete the argument. 

We will translate the claim into the symbolic calculus. 
Introduce symbolic letters $a,b$, and write 
\[ A = a_\ux^{rd-4p}, \quad B = b_\ux^d, 
\quad (\ux \, \uy) = x_0 \, y_1 - x_1 \, y_0. \] 
The rule for calculating transvectants 
symbolically is given in \cite[\S 49]{GrYo}; we use it to trace the steps 
from $T_1$ through $T_6$ for the calculation of  $u_r^{p,p'}(A \otimes B)$. 
Once this is done, we have the following statement to prove: 

\smallskip 

\begin{Lemma} \sl 
Given $r \ge 2$ and $0 \le p' \le \frac{(r+1)e}{2}$, there exists 
an integer $p$ in the range $0 \le p \le \frac{re}{2}$, such that 
the algebraic expression 
\[ \{ 
\Omega^{2p'} \; (\ux \, \uy)^{2p} \, a_{\ux}^{re-2p} \, a_{\uy}^{re-2p} \, 
b_{\ux}^{\,e} \, b_{\uy}^{\,e} \, \}|_{\uy:=\ux}
\] 
is not identically zero. 
\label{lemma.B} \end{Lemma}
\demo See Section~\ref{Feynman}. \qed 

\medskip 

At this point, modulo Lemmata~\ref{lemma.A} and~\ref{lemma.B}, 
the proofs of Theorems~\ref{surj.alpha_r} and~\ref{main.theorem} 
are complete. 
%%%%%%%%%%%%%%%%%%%%%%%%%%%%%%%%%%%%%%%%%%%%
\section{The combinatorics of Feynman diagrams} \label{Feynman}
We have kept the following presentation semi-formal, in order to avoid 
making the treatment cumbersome. Notwithstanding this, it is 
entirely rigorous as it stands. 
The reader looking for a stricly formal exposition of Feynman diagrams
should consult~\cite{Abdesselam} or~\cite{Fiorenza}, which 
implement Andr{\'e} Joyal's category-theoretic framework for combinatorial 
enumeration (see \cite{Bergeron,Joyal}). For our immediate purpose, 
let us simply say that a Feynman diagram is the combinatorial data 
needed to encode a complex tensorial expression built from a predefined collection
of elementary tensors, exclusively using contractions of
tensor indices. The word `tensor' is used here in the sense 
of a multidimensional analogue of a matrix, rather than the corresponding 
coordinate-free object from multilinear algebra.
Coordinates are needed in order to state the necessary definitions, 
but are almost never actually used in the computations.

\subsection{Diagrams and Amplitudes} \label{diag.ampl} 
Define the tensors 
\begin{equation} \ux=\left( 
\begin{array}{c} x_0 \\ x_1 \end{array} \right), \quad 
\uy=\left( \begin{array}{c} y_0\\y_1  \end{array} \right) 
\end{equation} 
made of formal indeterminates. Define the antisymmetric tensor 
$ \epsilon = \left( \begin{array}{rr} 0 & 1 \\ -1 & 0 \end{array} \right)$, 
and the symmetric tensor $Q$ which corresponds to the 
quadratic form $Q(\ux)=\ux^{\rm T} \, Q \, \ux$. 
Introduce the vectors of differential operators
\begin{equation}
\partial_\ux=\left( \begin{array}{c} \frac{\partial}{\partial x_0} \\ 
\frac{\partial}{\partial x_1} \end{array} \right), \quad  
\partial_\uy=\left( \begin{array}{c} \frac{\partial}{\partial y_0} \\
\frac{\partial}{\partial y_1} \end{array} \right) \end{equation}
We will use the following graphical notation for the entries of these tensors: 
\[
\figplace{dessin1}{0 in}{0 in}
\figplace{dessin2}{0 in}{0 in}
\figplace{dessin3}{0 in}{0 in}
\]
\begin{equation}
\figplace{dessin4}{0 in}{0 in}
\figplace{dessin5}{0 in}{0 in}
\figplace{dessin6}{0 in}{0 in}
\label{constituents}
\end{equation}
(The indices $\alpha,\beta$ belong to the set $\{0,1\}$).
We will obtain a `diagram' by assembling any number of these 
elementary pieces by gluing pairs of index-bearing lines; associated to it 
is an expression called the `amplitude' of the diagram.
Its rule of formation is as follows: 
introduce an index in $\{0,1\}$ for each glued pair of lines, 
take the product of the tensor entries corresponding to the 
different constituents from (\ref{constituents}) which appear in the 
diagram, and finally sum over all possible values of the indices. 
For instance, to the diagram 
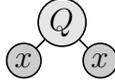
\[ \figplace{dessin7}{0 in}{-0.28 in} \] 
corresponds the amplitude 
$\sum_{\alpha,\beta \in \{0,1\}} x_\alpha \, Q_{\alpha \, \beta} \, x_\beta = Q(\ux)$, 
which is the quadratic form itself. Similarly, to 
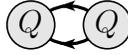
\[ \figplace{dessin8}{0 in}{-0.2 in} \] 
corresponds 
\[ \sum_{\alpha,\beta,\gamma,\delta \in \{0,1\}} 
Q_{\alpha \, \beta} \, \epsilon _{\alpha \, \gamma} 
\, \epsilon_{\beta \, \delta} \, Q_{\gamma \, \delta} 
  = 2 \, (Q_{00} \, Q_{11}-Q_{01}^2 ) \\
  = 2 \, \det(Q). \] 
Henceforth, whenever we write a diagram in an expression, it is the amplitude 
that is meant. Now the term $Q_{\alpha \, \beta}$ has an {\sl inner
structure}, related the notion of combinatorial plethysm 
(see~\cite{Bergeron,Joyal}).
Indeed, we can factor $Q$ as $Q(\ux)=R_1(\ux) \, R_2(\ux)$, where 
\begin{equation}
R_1=\left( \begin{array}{c} 
R_{1,0} \\ R_{1,1} \end{array} \right), \; 
R_2=\left( \begin{array}{c} R_{2,0} \\ R_{2,1} \end{array} \right)
\in \complex^2 \end{equation}
are dual to the homogeneous roots of $Q$. For any indices 
$\alpha$ and $\beta$, 
\[ Q_{\alpha \, \beta} = \frac{1}{2} \, 
\frac{\partial^2}{\partial x_\alpha \, \partial x_\beta} \, Q(\ux) 
= \frac{1}{2} \left( R_{1,\alpha} \, R_{2,\beta} + 
R_{1,\beta} \, R_{2,\alpha}  \right). \] 
We will write this more suggestively as
\begin{equation}
\figplace{dessin3}{0 in}{-0.25 in} = \figplace{dessin9}{0 in}{-0.25 in}
=\frac{1}{2} \figplace{dessin10}{0 in}{-0.5 in}
+\frac{1}{2} \figplace{dessin11}{0 in}{-0.5 in}
\label{Decomp} \end{equation}
This implies that 
\[ \figplace{dessin8}{0 in}{-0.2 in}
= \frac{1}{4} {{ \figplace{dessin12}{0 in}{-0.1 in}
}\atop{ \figplace{dessin13}{0 in}{-0.1 in} }} + \frac{1}{4}
{{ \figplace{dessin13}{0 in}{-0.1 in} }\atop{
\figplace{dessin12}{0 in}{-0.1 in} }} \] 
(Recall that reversing the direction of an $\epsilon$ arrow
introduces a minus sign, and therefore
\[ \figplace{dessin12bis}{0 in}{-0.1 in} =
\figplace{dessin13bis}{0 in}{-0.1 in} =0. \, ) \] 
Consequently, 
\[ \figplace{dessin8}{0 in}{-0.2 in} = -\frac{1}{2} \, \Delta, \] 
where 
\[ \Delta = \left(\figplace{dessin12}{0 in}{-0.1 in}\right)^2 \] 
is the discriminant of $Q$. 

\subsection{First Proof of Lemma \ref{lemma.A}}\label{first.proof.A}
Now write 
\[ F(\ux,\uy) = \Omega^{2p} \; Q(\ux)^e \, Q(\uy)^e, \] 
then $F(\ux,\ux) = F(\ux,\uy)|_{\uy = \ux}$ is the quantity we are 
interested in. Diagrammatically, $F(\ux,\uy)$ is equal to 
\begin{equation} 
\left(\figplace{dessin14}{0 in}{-0.09 in} \right)^{2p} 
\left(\figplace{dessin7}{0 in}{-0.26 in} \right)^e
\left( \figplace{dessin15}{0 in}{-0.26 in} \right)^e
\label{Fofxy} \end{equation} 
This is rewritten in terms of Feynman diagrams by summing over all ways to perform 
`Wick contractions' between $\frac{\partial}{\partial x}, \frac{\partial}{\partial y}$ 
on the one hand, and $x,y$ on the other hand (see e.g.~\cite{Abdesselam}). 
Once we let $\uy=\ux$, this condenses into the following sum over vertex-labelled 
bipartite multigraphs:
\begin{equation}
F(\ux,\ux) = \sum_G \, w_G \, \cA_G.  \label{SumG} 
\end{equation}
This is to be read as follows: we let $L$ and $R$ to be fixed sets of 
cardinality $e$ which label the $Q(\ux)$ and $Q(\uy)$ factors in (\ref{Fofxy})
respectively. Then a multigraph $G$ is identified with a matrix $(m_{ij})$ 
in $\NN^{L\times R}$. The quantity $w_G$ is the combinatorial weight and
$\cA_G$ is the amplitude of the Feynman diagram encoded by $G$.
Each $G$ entering into the sum satisfies the follwing conditions: 
\begin{itemize} 
\item 
$\sum\limits_{i\in L, j\in R} m_{ij}=2p$, 
\item For all $i \in L$, the number 
$l_i = \sum\limits_{j\in R} \, m_{ij}$ is $\le 2$, and 
\item For all $j \in R$, the number 
$c_j = \sum\limits_{i\in L} \, m_{ij}$ is $\le 2$. 
\end{itemize} 
The combinatorial weight is seen to be
\[ w_G= \frac{(2p)! \; 2^{2e}}
{\prod\limits_{i,j} (m_{ij})! \times \prod\limits_i (2-l_i)! \times
\prod\limits_j (2-c_j)!} \] 
The amplitude $\cA_G$ factors over the connected components
of $G$. These components are of four possible types:
cycles containing an even number of $\epsilon$ arrows
of alternating direction, chains with both endpoints in $L$,
chains with both endpoints in $R$, and finally chains with one 
endpoint in $L$ and another in $R$. However, the contribution from the 
last type is zero. Indeed, such a chain contains an odd 
number of $\epsilon$ arrows, and therefore its amplitude changes sign
if we reverse the orientations on all the arrows. 
But the last operation, followed by a rotation of $180^\circ$, puts
the chain back in its original form. For instance, 
\medskip 

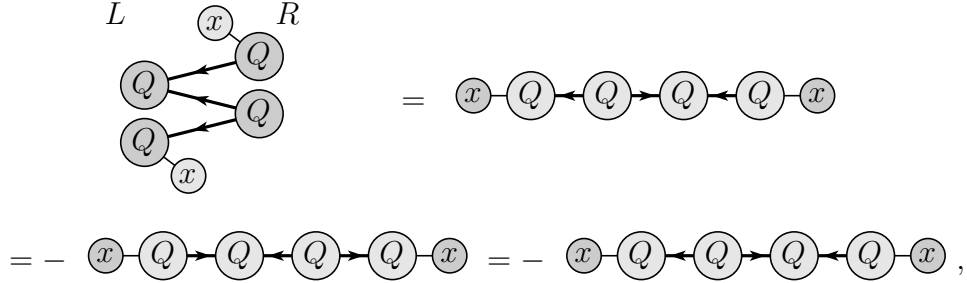
\[ \begin{aligned} 
{} & \figplace{dessin16}{0 in}{-0.45 in} 
= \figplace{dessin17}{0 in}{-0.18 in}  \\ 
& = -\figplace{dessin18}{0 in}{-0.18 in} 
= -\figplace{dessin17}{0 in}{-0.18 in}, 
\end{aligned} \] 
and hence this expression vanishes. 
Now we can use the inner structure of $Q$ to 
calculate the other three amplitudes. 
Given a cycle of even length $2m$, we incorporate the decomposition (\ref{Decomp})
at each $Q$ vertex. This produces a sum of $2^{2m}$ terms, 
all but two of which vanish. 
Indeed, suppose we have chosen the precise connections between
the `inner' and `outer' part of what was a particular $Q$ vertex. 
Then, since the vanishing factors 
\[ \figplace{dessin12bis}{0 in}{-0.09 in}, 
\text{and} \figplace{dessin13bis}{0 in}{-0.09 in}
\] 
are to be avoided, the connections for the remaining vertices are forced. 
Moreover, the alternating pattern for the orientations of the 
$\epsilon$ arrows implies that we collect an equal number $m$ of
either of the following factors: 
\[ \figplace{dessin12}{0 in}{-0.09 in} \quad 
\text{and} \figplace{dessin13}{0 in}{-0.09 in}
\] 
As a result, the amplitude of the cycle is exactly 
$(-\Delta)^m \, 2^{1-2m}$.
Similarly, a chain with both endpoints in $L$ (or both in $R$) 
and with a necessarily even number $2m$ of $\epsilon$ 
arrows (and thus $2m+1$ of $Q$ vertices) gives an amplitude
\[ \frac{2}{2^{2m+1}} \times 
\figplace{dessin19}{0 in}{-0.09 in}
\left( \figplace{dessin13}{0 in}{-0.09 in} 
\figplace{dessin12}{0 in}{-0.09 in} \right)^m 
\figplace{dessin20}{0 in}{-0.09 in} \]
\[ = 2^{-2m} (-\De)^m Q(\ux). \] 
Therefore, an easy count shows that the amplitude of a bipartite 
multigraph $G$ in (\ref{SumG}) is
\[ \cA_G = 2^{\, \cC(G)-2p} \, Q(\ux)^{2e-2p} \, (-\Delta)^p, \] 
where $\cC(G)$ is the number of cycles in $G$. Finally,
\[ F(\ux,\ux)=\cN_{e,p}^{\, \rm I} \, Q(\ux)^{2e-2p} \, (-\Delta)^p,
\] 
where $\cN_{e,p}^{\, \rm I}$ denotes the sum 
\begin{equation}
\sum_{G} \; \frac{(2p)! \times 2^{\, 2e-2p+\cC(G)}} 
{\prod_{i,j} (m_{ij})! \times \prod_i (2-l_i)! \times \prod_j (2-c_j)!} 
\label{NIdef} \end{equation}
The sum is quantified over all $G=(m_{ij})$ satisfying the 
three constraints above, and the additional constraint 
that there is no connected connected component which is a chain starting
in $R$ and ending in $L$. 
It is not difficult to see that given $e\ge 1$ and $0\le p\le e$,
there always exists such a graph $G$. For instance, take $G$ corresponding
to a matrix having $p$ of its diagonal entries set equal to $2$ and
zeroes elsewhere. Hence, 
$\cN_{e,p}^{\, \rm I}>0$ which proves Lemma \ref{lemma.A}. \qed
\begin{Remark} \rm 
The factor of $+2$ per cycle in (\ref{NIdef}) should be 
contrasted with the $-2$ factor in Penrose's original definition 
of spin networks~\cite{Penrose}. This intuitively suggests 
that Penrose's construction might be a Fermionic or `negative dimensional'
analogue of covariants of binary forms. 
\end{Remark}

\subsection{Second proof of Lemma \ref{lemma.A}}\label{second.proof.A}
Let $p,q,k$ be nonnegative integers, with
$k\le 2 \min\{ p,q  \}$. Let $Q$ be a binary quadratic with 
discriminant $\Delta$ (normalized as in the previous section). 
We will calculate the transvectant $\T = (Q^p,Q^q)_k$ precisely. 
The special case $p=q=e$ gives another proof of Lemma~\ref{lemma.A}. 
\begin{Proposition}\label{prop.hyper} \sl 
If $k$ is odd, then $\T =0$. If $k=2m$ is even, then
\begin{equation}
\T = Q^{p+q-2m} \, (-\Delta)^m \times \cN^{\, \rm II}_{p,q,m}, 
\label{Lemmadixon1} \end{equation} 
where 
\begin{equation} \cN^{\, \rm II}_{p,q,m} = 
\frac{p! \, q! \, (2m)! \, (p+q-m)! \, (2p-2m)! \, (2q-2m)!} 
{(2p)! \, (2q)! \, m! \, (p+q-2m)! \, (p-m)! \, (q-m)!}. \label{Lemmadixon2} 
\end{equation} \end{Proposition} 
\demo
We specialize the quadratic form to $Q(\ux)=x_0 \, x_1$, for which
$\De=1$. Since $(Q^q,Q^p)_k=(-1)^k (Q^p,Q^q)_k$, we may assume 
$p\le q$. Now expand $\Omega^k$ by the binomial theorem. 
By definition (\ref{trans.defn}), 
\begin{equation}
\begin{aligned} {} & (Q^p,Q^q)_k= \frac{(2p-k)! \, (2q-k)!}{(2p)! \, (2q)!} \, \times \\ 
& \left. \sum_{i=0}^k (-1)^i \binom{k}{i}  
{\lp\frac{\partial}{\partial x_0}\rp}^{k-i}
{\lp\frac{\partial}{\partial y_1}\rp}^{k-i}
{\lp\frac{\partial}{\partial x_1}\rp}^{i}
{\lp\frac{\partial}{\partial y_0}\rp}^{i}
\times x_0^p \, x_1^p \, y_0^q \, y_1^q \right|_{\uy:=\ux \, .} 
\end{aligned} \end{equation}
After differentiating and letting $\uy:=\ux$, this reduces to 
\[ (Q^p,Q^q)_k=C_{p,q,k}\times W_{p,q,k}, \] 
where 
\[ C_{p,q,k} = \frac{(2p-k)! \, (2q-k)! \, k! \, (p!)^2 \, (q!)^2} 
{(2p)! \, (2q)!} \times x_0^{p+q-k} \, x_1^{p+q-k},  \] 
and 
\begin{equation}
W_{p,q,k} = \sum_{\max \{0,k-p\}}^{\min\{k,p\}} 
\frac{(-1)^i}{i!(k-i)!(p-i)!(q-i)!(p-k+i)!(q-k+i)!}
\label{Wsum}. \end{equation}
Up to a numerical factor, (\ref{Wsum}) is Van der Waerden's
formula for Wigner's $3j$-symbols (see~\cite{Biedenharn1}). 
We now have two cases to consider.

\medskip 

\noindent{\bf First case :} Assume $0\le k\le p$. 
Using Pochammer's symbol
\[ (a)_i := a(a+1)\cdots (a+i-1), \] 
and the obvious identities $(a+i)!=a! (a+1)_i$ and
$(a-i)!=\frac{(-1)^i a!}{(-a)_i}$, we can write
\[ W_{p,q,k}=\frac{1}{k! \, p! \, q! \, (p-k)! \, (q-k)!} \, 
\sum_{i=0}^k \, \frac{(-k)_i(-p)_i(-q)_i}{i!(p-k+1)_i(q-k+1)_i}, 
\] or
\begin{equation}
W_{p,q,k}=\frac{1}{k! \, p! \, q! \, (p-k)! \, (q-k)!} \; 
{ }_3 \bF_2 \left[ \begin{array}{c} -k,-p,-q \\ 
p-k+1, q-k+1 \end{array} ; 1 
\right] \label{doublev}. \end{equation} 
The hypergeometric series appearing in this formula can be evaluated by 
Dixon's summation theorem (see~\cite[p.~52]{Slater}). It gives the formula 
\[ \begin{aligned} {}_3 \bF_2 & \left[
\begin{array}{c} a,b,c \\ 1+a-b, 1+a-c \end{array} ; 1 
\right] = \\ 
& \frac{\Gamma(1+\frac{1}{2}a) \, \Gamma(1+\frac{1}{2}a-b-c) \, 
\Gamma(1+a-b) \, \Gamma(1+a-c)} {\Gamma(1+a) \, \Gamma(1+a-b-c) \, 
\Gamma(1+\frac{1}{2}a-b) \, \Gamma(1+\frac{1}{2}a-c)}, 
\end{aligned} \] 
which is valid in the domain of analyticity $\Re({1+\frac{1}{2}a-b-c})>0$.
We would like to choose $a=-k$, $b=-p$ and $c=-q$, hence 
we rewrite the factor $\frac{\Gamma(1+\frac{1}{2}a)}{\Ga(1+a)}$ 
as 
\begin{equation}
\frac{\pi}{\Ga(-\frac{a}{2})\sin(-\frac{\pi a}{2})}
\times\frac{\Ga(-a)\sin(-\pi a)}{\pi}=
\cos(\frac{\pi a}{2}) \, \frac{\Ga(-a+1)}{\Ga(-\frac{a}{2}+1)}
\end{equation}
before specializing $a,b,c$.
That is, we use Dixon's theorem in the form 
\begin{equation} \begin{aligned} { }_3 \bF_2
& \left[ \begin{array}{c} a,b,c \\ 1+a-b, 1+a-c \end{array}
; 1 \right]=\cos \left( \frac{\pi a}{2}\right) \times \\ 
& \frac{\Gamma(1-a) \, \Gamma(1+\frac{1}{2}a-b-c) \, 
\Gamma(1+a-b) \, \Gamma(1+a-c)} {\Gamma(1-\frac{a}{2}) \, \Gamma(1+a-b-c) \, 
\Gamma(1+\frac{1}{2}a-b) \, \Gamma(1+\frac{1}{2}a-c)}
\label{Dixon} \end{aligned} \end{equation}
Now let $a=-k$, $b=-p$ and $c=-q$. Then, since 
$0\le k\le p\le q$, all the arguments of the Gamma function are 
strictly positive. 

If $k$ is odd, the cosine factor vanishes, and hence 
so does $\T$. (This vanishing has a different explanation in 
the context of the first proof above: since there is an 
odd number of arrows, there must exist a chain joining $L$ to $R$.) 
If $k=2m$ is even, then 
formulae (\ref{doublev}) and (\ref{Dixon}) imply that 
\begin{equation}
W_{p,q,2m}=(-1)^m
\frac{(p+q-m)!}{p! \, q! \, m! \, (p+q-2m)! \, (p-m)! \, (q-m)!} 
\end{equation}
which implies (\ref{Lemmadixon1}) for 
the quadratic form $Q(\ux)=x_0 \, x_1$. Since a generic quadratic form 
lies in the $GL_2(\complex)$ orbit of $x_0 \, x_1$, the formula is 
proved in general. 

\smallskip 

\noindent{\bf Second case :} Assume $k>p$. We make a change of 
index $i=k-p+j$, then (\ref{Wsum}) becomes
\[ \begin{aligned} 
{} & W_{p,q,k} =\sum_{j=0}^{2p-k}\\ 
& \frac{(-1)^{k-p+j}} 
{j! \, (p-j)! \, (2p-k-j)! \, (p+q-k-j)! \, (k-p+j)! \, (q-p+j)!}
\end{aligned} \] 
Once again, this can be rewriten as an  ${ }_3 \bF_2$ hypergeometric 
series to which Dixon's theorem applies. 
\[ \begin{aligned} {} &
W_{p,q,k}=\frac{(-1)^{k+p}}{p! \, (2p-k)!\,(p+q-k)!\,(k-p)!\,(q-p)!}\\
& \times\ { }_3 \bF_2 \left[
\begin{array}{c} -2p+k, -p, -p-q+k \\ k-p+1, q-p+1 \end{array} ; 1 
\right]. \end{aligned} \] 
Now we apply Dixon's theorem in the modified form (\ref{Dixon}), 
with $a=-2p+k$, $b=-p$, $c=-p-q+k$ and conclude as before. The 
proposition (and Lemma~\ref{lemma.A}) are proved. \qed

We obtain a closed formula for the weighted graph enumeration (\ref{NIdef}) 
by comparing both proofs of the Lemma: 
\[ \cN_{e,p}^{\, \rm I} =\frac{(2e)!^2}{(2e-2p)!^2} \; 
\cN_{e,e,p}^{\, \rm II}. \] 

\smallskip 

\subsection{Proof of Lemma \ref{lemma.B}}
Let $r$, $e$, $p'$ and $p$ be integers satisfying
$r\ge 2$, $e\ge 1$, $0\le 2p'\le (r+1)e$ and $0\le 2p\le re$. Let 
\[ a=\left( \begin{array}{c} a_0 \\ a_1 \end{array} \right), \quad 
b=\left( \begin{array}{c} b_0 \\ b_1 \end{array} \right) 
\] 
be two elements of $\complex^2$ and
\[ \ux=\left( \begin{array}{c} 
x_0 \\ x_1 \end{array} \right), \quad 
\uy=\left( \begin{array}{c} y_0 \\ y_1 \end{array} \right) 
\] be two vectors of indeterminates. The quantity we would like to 
calculate is
\[ G(\ux) =  
\{ \left. \Omega^{2p'} \; (\ux \, \uy)^{2p} \, 
a_\ux^{re-2p} \, a_\uy^{re-2p} \, b_\ux^e \, b_\uy^e \}
\right|_{\uy:=\ux \,.}  \] 
or, in matrix notation,
\[ G(\ux)= \left. \{
[\partial_\ux^{\rm T} \, \epsilon \, \partial_\uy]^{2p'}
(\ux^{\rm T} \, \epsilon \, \uy)^{2p} \, 
(\ux^{\rm T} \, a \, a^{\rm T} \, \uy)^{re-2p} \, 
(\ux^{\rm T} \, b \, b^{\rm T} \, \uy)^{e} \} 
\right|_{\uy:=\ux \,.} \] 
Introduce two new vectors of auxiliary variables
\[ \phi=\left( \begin{array}{c} 
\phi_0 \\ \phi_1 \end{array} \right), \quad 
{\overline\phi}=\left( \begin{array}{c} 
{\overline\phi}_0\\ {\overline\phi}_1 \end{array}
\right) \] 
and rewrite $G(\ux)$ as 
\[ \begin{aligned} 
{} & (2p')! \, (2p)! \, (re-2p)! \, e! \, 
\frac{[\partial_\ph^{\rm T} \, \epsilon \, \partial_{\overline\ph}]^{2p'}}{(2p')!}
\frac{[(\phi+\ux)^{\rm T} \, \epsilon \, ({\overline \phi}+\ux)]^{2p}}{(2p)!} \\
& \times 
\frac{[({\overline\phi}+\ux)^{\rm T} \, a \, a^{\rm T}(\phi+\ux)]^{re-2p}}
{(re-2p)!} \times \left. 
\frac{[({\overline\phi}+\ux)^{\rm T} \, b \, b^{\rm T}(\phi+\ux)]^{e}} {e!}
\right|_{\phi,\overline\phi:=0 \,.}
\end{aligned} \] 
\subsection{A `Gaussian integral' on $\complex^2$} 
We now introduce a term $\cZ$, which can be seen as the 
combinatorial algebraic avatar of a Gaussian integral on $\complex^2$. 
(Compare~\cite{Abdesselam,Bargmann}, where the ${\overline \phi}$ 
are actual complex conjugates of the $\phi$.)

We will write $\complex \, [[\xi_1,\xi_2,\dots]]$ for 
the ring of formal power series in variables $\xi_1,\xi_2$ etc. Define 
$M = v \, a \, a^{\rm T}+ w \, b \, b^{\rm T}$, 
a $2 \times 2$ matrix over $\complex[[v,w]]$. Let 
\[ S = ({\overline\phi}+\ux)^{\rm T}(-u \, \epsilon +M)(\phi+\ux)  \in 
\complex[[\phi_0, \phi_1,{\overline\phi}_0,{\overline\phi}_1,x_0,x_1,h,u,v,w]] 
\] 
and define 
\[ \cZ = \left. \left\{ 
\sum_{n\ge 0} \, \frac{h^n}{n!} \, 
[\partial_\phi^{\rm T} \, \epsilon \, \partial_{\overline\phi}]^{n} \, 
e^S \right\} \right|_{\phi,\overline\phi:=0}  \in 
\complex \, [[x_0,x_1,h,u,v,w]]. \] 
Then we have \[ G(\ux)=(2p')! \, (2p)! \, (re-2p)! \, e! \; 
[h^{2p'} \, u^{2p} \, v^{re-2p} \, w^e]_\cZ, \] 
where $[h^{2p'} \, u^{2p} \, v^{re-2p} \, w^e]_\cZ \in\complex \, [[x_0,x_1]]$
denotes the coefficient of the monomial $h^{2p'} \, u^{2p} \, v^{re-2p} \, w^e$ in
$\cZ$. With obvious notations, one can rewrite $\cZ$ as
\[ \cZ= \left.
\exp(h \, \partial_\phi^{\rm T} \, \epsilon \, \partial_{\overline\phi}) \, 
\exp({\overline\phi}^{\rm T} \, A \, \phi+J^{\rm T} \phi+{\overline\phi}^{\rm T} \, K
+S_0) \right|_{\phi, \overline\phi:=0}, \] 
with 
\[ \begin{array}{ll} 
A = -u \, \epsilon +M, & J^{\rm T} = -u \, \ux^{\rm T} \, \epsilon+\ux^{\rm T} \, M, \\ 
K=-u \, \epsilon \, \ux+M \, \ux, & 
S_0 = v \, (\ux^{\rm T}\, a)^2+ w \, (\ux^{\rm T} \, b)^2. \end{array} \] 
Therefore $\cZ=e^{S_0} \, {\widetilde\cZ}$ with
\[ {\widetilde\cZ} = \left.
\exp(h \, \partial_\phi^{\rm T} \, \epsilon \, \partial_{\overline\phi}) \, 
\exp({\overline\phi}^{\rm T} \, A \, \phi+J^{\rm T}\, \phi+{\overline\phi}^{\rm T} \, K)
\right|_{\phi, \overline\phi:=0 \,.} \] 

Now ${\widetilde\cZ}$ can be expressed as a sum over Feynman diagrams, 
built as in Section \ref{diag.ampl}, from the following pieces 
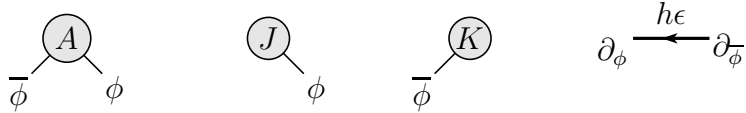
\[ \figplace{dessin21}{0 in}{0 in}
\figplace{dessin22}{0 in}{0 in} \figplace{dessin23}{0 in}{0 in}
\figplace{dessin24}{0 in}{0 in} \]
by plugging the $\partial_\phi$ onto the $\phi$, and  
the $\partial_{\overline\phi}$ onto the $\overline\phi$ in all possible ways.

More precisely, given any finite set $E$, we define a Feynman diagram
on $E$ as a sextuple $\cF=(E_\phi,E_{\overline\phi}, \pi_A,\pi_J,\pi_K,\cC)$, 
where \begin{itemize} 
\item $E_\ph$, $E_{\Br\ph}$ are subsets of $E$, 
\item $\pi_A$, $\pi_J$, $\pi_K$ are sets of subsets of $E$, and 
\item $\cC$ is a map $E_{\overline\phi} \longrightarrow E_\phi$; 
\end{itemize} 
satisfying the following axioms:
\begin{itemize}
\item
$E_\ph$ and $E_{\Br\ph}$
have equal cardinality and they form a two set partition of $E$.
\item
The union of the elements in $\pi_A$, that of elements in $\pi_J$, 
and that of elements in $\pi_K$ form a three set partition of $E$.
\item 
$\cC$ is bijective.
\item
Every element of $\pi_A$ has two elements, one in $E_\phi$ and 
one in $E_{\overline\phi}$.
\item
Every element of $\pi_J$ has only one element which lies in $E_\phi$.
\item
Every element of $\pi_K$ has only one element which lies 
in $E_{\overline\phi}$.
\end{itemize}
The set of Feynman diagrams on $E$
is denoted by ${\mathsf{Fey}}(E)$. Given a Feynman 
diagram $\cF$ on $E$ and a bijective map $\si:E\rightarrow E'$, there 
is a natural way to transport $\cF$ along $\si$ in order to 
obtain a Feynman diagram $\cF'={\mathsf{Fey}}(\si)(\cF)$ on $E'$. 
Hence $E\rightarrow {\mathsf{Fey}}(E)$
defines an endofunctor of the groupoid category of finite
sets with bijections (cf.~\cite{Abdesselam,Bergeron,Fiorenza,Joyal}). 

\begin{Example} \rm 
Let $E=\{1,2,\ldots,8\}$, $E_\ph=\{1,2,3,4\}$, 
$E_{\Br\ph}=\{5,6,7,8\}$, $\pi_A=\{ \{2,6\},\{3,7\}, \{4,8\} \}$, 
$\pi_J=\{\{1\}\}$, $\pi_K=\{\{5\}\}$, and $\cC$ given by $\cC(5)=1$, $\cC(6)=3$,
$\cC(7)=4$ and $\cC(8)=2$. This corresponds to the diagram
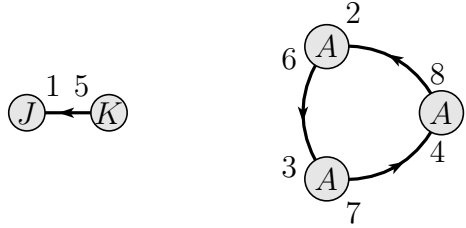
\[
\figplace{dessin25}{0 in}{0.8 in}
\figplace{dessin26}{0 in}{0 in}
\]
where we put the elements of $E$ next to the corresponding half-line.
The amplitude of such a pair $(E,\cF)$ is 
\[ 
\cA(E,\cF)=(J^{\rm T} \, (h \, \epsilon) \, K) 
\times \text{trace}([h \, \epsilon \, A]^3). \] 
\end{Example} 
There is a natural equivalence relation between pairs
of finite sets equiped with a Feynman diagram. It is given by letting
$(E,\cF)\sim(E',\cF')$ if and only if there exists a bijection
$\si:E\rightarrow E'$ such that $\cF'={\mathsf{Fey}}(\si)(\cF)$.
The automorphism group $\text{Aut}(E,\cF)$
of a pair $(E,\cF)$ is the set of bijections $\si:E \longrightarrow E$ 
such that ${\mathsf{Fey}}(\si)(\cF)=\cF$. Now, 
\[ {\widetilde\cZ}=\sum_{[E,\cF]} \, \frac{\cA(E,\cF)}{|\text{Aut}(E,\cF)|}
\]
where the sum is quantified over equivalence classes of 
pairs $(E,\cF)$. The term $\cA(E,\cF)$ is the amplitude, 
and $|\text{Aut}(E,\cF)|$ is the cardinality of the automorphism group. 
We leave it to the reader to check 
(otherwise see~\cite{Abdesselam,Fiorenza}) that
\[ \begin{aligned} 
\log {\widetilde\cZ} & =  \sum_{[E,\cF]\ \rm connected} \,  
\frac{\cA(E,\cF)}{|\text{Aut}(E,\cF)|} \\ 
& = \sum_{n\ge 1} \, \frac{1}{n} \, \text{trace}((h \, \epsilon \, A)^n) +
\sum_{n\ge 0} \, J^{\rm T} \, (h \, \epsilon \, A)^n \, (h \, \epsilon) \, K. 
\end{aligned} \] 
(This uses the fact that the only connected diagrams are pure $A$-cycles 
or $A$-chains joining a $J$ to a $K$ vertex.) Hence 
\[ {\widetilde\cZ}=\frac{1}{\det(I-h \, \epsilon \, A)} 
\exp(J^{\rm T}\, (I-h\ep A)^{-1} \, (h \, \epsilon) \, K). 
\] 
After straightforward but tedious computations with $2 \times 2$ 
matrices (which we spare the reader), one gets 
\[ \cZ=\frac{1}{(1-h \, u)^2+h^2 \, v \, w \, (a^{\rm T} \, \epsilon \, b)^2}
\exp \left( 
\frac{v \, (\ux^{\rm T} \, a)^2+w \, (\ux^{\rm T} \, b)^2}
{(1-h \, u)^2+h^2 \, v \, w \, (a^{\rm T} \, \epsilon \, b)^2} \right), 
\] 
or in classical notation
\[ \cZ=\frac{1}{(1-h\, u)^2+h^2 \, v \, w \, (a \, b)^2}
\exp\left(\frac{v \, a_\ux^2+w \, b_\ux^2}{(1-h\,u)^2+h^2 \, v \, w \, (a \, b)^2}
\right). \] 
Expanding this, 
\[ \begin{aligned} 
\cZ & =\sum_{\mu\ge 0} \, \frac{1}{\mu!} (v \, a_\ux^2+w \, b_\ux^2)^{\mu} 
\left( 1-2 \, h \, u+h^2 \, u^2+h^2 \, v \, w \, (a \, b)^2 \right)^{-(\mu+1)} \\ 
& = \sum_{\mu,\nu\ge 0} \frac{(-1)^\nu (\mu+\nu)!}{\mu!^2 \, \nu!}
(v \, a_\ux^2+w \, b_\ux^2)^{\mu} 
(-2 \, h \, u+h^2 \, u^2+h^2 \, v \, w \, (a \, b)^2)^{\nu} \\ 
& =\sum_{\stackrel{m,n}{\alpha,\beta,\gamma}\ge 0}
\frac{(-1)^{\alpha+\beta+\gamma}(m+n+\al+\beta+\ga)!} 
{(m+n)! \, m! \, n! \, \alpha! \, \beta! \, \gamma!} \; \times \\ 
& \qquad (v \, a_\ux^2)^m (w \, b_\ux^2)^n \, (-2h \, u)^\alpha (h^2 \, u^2)^\beta 
\, (h^2 \, v \, w \, (a \, b)^2)^\gamma \\ 
& =\sum_{\stackrel{m,n}{\alpha,\beta,\gamma} \ge 0}
\frac{(-1)^{\beta+\ga} 2^\al (m+n+\al+\beta+\ga)! }
{(m+n)! \, m! \, n! \alpha! \, \beta! \, \gamma!} \; \times \\ 
& \qquad 
h^{\alpha+2 \, \beta+2 \, \gamma} \, 
u^{\alpha+2 \, \beta} \, v^{m+\gamma} \, w^{n+\gamma} \, a_\ux^{2m} \, 
b_\ux^{2n} \, (a \, b)^{2 \, \gamma}. 
\end{aligned} \] 
The coefficient of $h^{2p'} \, u^{2p} \, v^{re-2p} \, w^e$
is a sum over the single index $\beta$, $0\le \beta\le p$,
as a result of solving for $\alpha=2p-2\beta$, $\gamma=p'-p$,
$m=re-p'-p$, $n=e-p'+p$. Therefore
\[ G(\ux)=\cN_{r,e,p',p}^{\, \rm III} \; 
a_\ux^{2(re-p'-p)} \, b_\ux^{2(e-p'+p)} \, (a \, b)^{2(p'-p)},
\] 
where
\[ \begin{aligned} {} & \cN_{r,e,p',p}^{\, \rm III} = 
\bbbone_{\left\{ p'-p\ge 0, \, e-p'+p\ge 0, \, re-p'-p\ge 0 \right\}} \; 
\times \\ 
& \frac{(-1)^{p'-p} \, (2p)! \, (2p')! \, (re-2p)! \, e!}
{(p'-p)! \, (e-p'+p)! \, (re-p'-p)! \, ((r+1)e-2p')!} \times \cJ_{s,p}. 
\end{aligned} \] 
Here $\bbbone_{\{\cdots\}}$ denotes the characteristic function
of the condition between braces, and
\[ \cJ_{s,p} = \sum_{\beta=0}^p \, 
\frac{(-1)^\beta \, 2^{2p-2 \, \beta} \, (s+2p-\beta)!}
{(2p-2 \, \beta)! \, \beta!} \] 
with $s = (r+1)e-p'-p$. Note that $s\ge e$ whenever 
the characteristic function is nonzero. Now $\cJ_{s,p}$ 
can be rewritten as a Gauss hypergeometric series 
and can be summed by the Chu-Vandermonde theorem (see~\cite[p.~28]{Slater}). 
The result is 
\[ \cJ_{s,p}=\frac{(s+p)! \, (s+\frac{3}{2})_p}{p! \, (\frac{1}{2})_p}  \] 
As a result, the characteristic function alone dictates whether 
$\cN_{r,e,p',p}^{\, \rm III}$ vanishes or not. 

Now for $r\ge 2$, $e\ge 1$ and $0\le p'\le \frac{(r+1)e}{2}$, 
it is easy to see that one can always find an integer $p$ 
with $0\le p\le \frac{re}{2}, p'-p\ge 0,e-p'+p\ge 0$ and $re-p'-p\ge 0$.
Indeed, take $p=p'$ if $0\le p'\le \frac{re}{2}$,
and otherwise take $p=p'-e$ if $\frac{re}{2}< p'\le \frac{(r+1)e}{2}$.
In either case, this ensures that $G(\ux)$ does not vanish identically, 
which proves Lemma~\ref{lemma.B}. \qed

\smallskip 

\noindent The proofs of Theorems~\ref{surj.alpha_r} 
and~\ref{main.theorem}  are now complete. 
\subsection{Proof of Proposition \ref{magic.squares}}
Write $Q = \prod\limits_{i=1}^r ( l_{i,0} \, x_0 + l_{i,1} \, x_1)$, where 
the $l_{i,-}$ are indeterminates. 
By Proposition~\ref{key.transvect}, the polynomial 
\[ (\prod_{i=1}^r  l_i^e,\prod_{j=1}^r l_j^e)_{2p} \]
is not identically zero. Now, as in section \ref{first.proof.A}, 
one can calculate the previous transvectant via the expression 
\[ \frac{(re-2p)!^2}{(re)!^2} \left. 
\Omega^{2p} (\prod_{i=1}^r  \, l_i(\ux)^e) 
(\prod_{j=1}^r \, l_j(\uy)^e) \right|_{\uy:=\ux} \]
by summing over the derivative actions. This
generates a sum over bipartite graphs between two sets
of $r$ elements which separately label the linear forms in each
of the two products. The valences of the vertices are bounded by
$e$ and the total number of edges is $2p$. Thus, 
\[ \begin{aligned}
{} & (\prod_{i=1}^r  l_i^e,\prod_{j=1}^r l_j^e)_{2p} 
=\frac{(re-2p)!^2}{(re)!^2}\times \\ 
& \sum_{N\in \cN} w_N \, (\prod_{1\le i,j\le r}
(l_i \, l_j)^{n_{ij}} \prod_{1\le i\le r} l_i(\ux)^{e-l(N)_i}
\prod_{1\le j\le r} l_j(\ux)^{e-c(N)_i}). \end{aligned} \]
Here $\cN$ is the set of $r\times r$ matrices
$N=(n_{ij})_{1\le i,j\le r}$ with nonnegative integer entries, 
such that
\begin{itemize} 
\item $ \sum\limits_{1\le i, j\le r} n_{ij}=2p$, 
\item for all $1\le i\le r$, the integer 
$l(N)_i = \sum\limits_{1\le j\le r} n_{ij}$ is $\le e$, 
\item 
for all $1\le j\le r$, the integer 
$c(N)_j = \sum\limits_{1\le i\le r} n_{ij}$ is $\le e$. 
\end{itemize} 
The combinatorial weight $w_N$ is given by
\[ 
\frac{(2p)!}{\prod\limits_{1\le i,j\le r} n_{ij}!}
\times \prod_{1\le i\le r} \frac{e!}{(e-l(N)_i)!}
\times \prod_{1\le j\le r} \frac{e!}{(e-c(N)_j)!}
\] and the bracket factors $(l_i \, l_j)$ stand for
$l_{i,0} \, l_{j,1}-l_{i,1} \, l_{j,0}$.
For given edge multiplicities recorded in the matrix $N$,
the combinatorial weight counts in how many ways one can obtain the 
correponding configuration by differentiating. Among the $2p$ lines, 
we have to choose which ones are assigned to each pair of vertices $(i,j)$, 
this gives the first multinomial factor. Then one has to specify the 
connections at each vertex, and this accounts for the other two factors.

Let us border the matrix $N$ by adding an extra row and column 
to make an $(r+1)\times(r+1)$ matrix $M$. Insert the `defect' numbers 
$\{e-c(N)_j\}$ into the $(r+1)$-th row, the $\{e-l(N)_i\}$ 
into the $(r+1)$-th column, and a $0$ in the bottom right corner. 
Then, keeping the notation of Proposition~\ref{magic.squares}, we have 
\[ \begin{aligned} 
{} & (\prod_{i=1}^r  \, l_i^e,\prod_{j=1}^r \, l_j^e)_{2p}
=\frac{(re-2p)!^2 \, (2p)! \, e!^{2r}}{(re)!^2} \, \times\\
& \sum_{M\in \cM}
\frac{\prod\limits_{1\le i,j\le r} (l_i \, l_j)^{m_{ij}}
\prod\limits_{1\le i\le r} l_i(\ux)^{m_{i,r+1}}
\prod\limits_{1\le j\le r} l_j(\ux)^{m_{r+1,j}}}
{\prod\limits_{1\le i,j\le r+1}m_{ij}!} 
\end{aligned} \] 
Dehomogenize the last expression by 
substituting $l_{i,0}=z_i$, $l_{i,1}=1$ for every $1\le i\le r$, and 
$x_0=-1, x_1=t$. This is a numerical multiple of the 
expression~(\ref{tau}), hence we have shown 
Proposition~\ref{key.transvect}. \qed 

%%%%%%%%%%%%%%%%%%%%%%%%%%%%%%%%%%%%%%%%%%%%%%
\section{Binary Forms} \label{binarycase}
Let $n=1$, then $X = X_{(1,d)}$ is the locus of 
degree $d$ binary forms which are $e$-th powers of quadratic forms. 
The following result for the $d=4$ case is classical 
(see \cite[\S 3.5.2]{Glenn}): 
\begin{Proposition} \sl 
A binary quartic $F$ lies in $X$, iff $((F,F)_2,F)_1 = 0$. 
\end{Proposition} 
We will generalize this to any (even) $d$ by identifying  
the covariants of binary $d$-ics which correspond to the generators 
of the ideal 
\[ I_X \subseteq \bigoplus\limits_{i \ge 0} S_i(S_d). \] 
Since $m_0=3$ and $(I_X)_2=0$, all the generators are in degree $3$. 
(It follows that the graded minimal resolution of $I_X$ is linear, 
however we will make no use of this.) 

\subsection{Cubic Covariants} 
We have a decomposition 
\[ (I_X)_3 = \bigoplus\limits_m \, (S_m \otimes \complex^{\nu_m}) 
\subseteq S_3(S_d) \] 
of irreducible $SL_2$-modules. 
Each equivariant inclusion $S_m \subseteq S_3(S_d)$ corresponds to 
a covariant of degree $3$ and order $m$ of binary $d$-ics. (This 
correspondence is explained in \cite{Chipalkatti1}.) Thus we have 
$\nu_m$ linearly independent cubic covariants of order $m$ which vanish 
on $X$. If we can list every covariant which occurs this way, 
then $I_X$ is completely specified. 

For the generic binary $d$-ic $F$, define 
\[ \U(i,j) = ((F,F)_{2i},F)_j, \] 
a covariant of degree $3$ and order $3d-4i-2j$. Unless the conditions 
\[ 0 \le i \le e, \quad 0 \le j \le \min \, \{d,2d-4i\} \qquad (\dagger), \] 
hold, $\U(i,j)$ is identically zero; hence we always assume that the 
pair $(i,j)$ satisfies $(\dagger)$. If $j$ is even, then write $j = 2k$, and 
define a rational number 
\[ \mu_{i,j} = (-1)^{i+k} \, \cN^{\, \rm II}_{e,e,i} \times 
\cN^{\, \rm II}_{2e-2i,e,k},  \] 
where $\cN^{\, \rm II}$ is defined by formula~(\ref{Lemmadixon2}).
If $(i,j),(\tilde i, \tilde j)$ are two pairs such that 
$2i + j = 2 \tilde i + \tilde j$, and 
$j,\tilde j$ are even, then define 
\[ \Phi(i,j,\tilde i, \tilde j) = \mu_{\tilde i, \tilde j} \;  \U(i,j) - 
\mu_{i,j} \;  \U(\tilde i,\tilde j). \] 
Now consider the following set of covariants: 
\begin{equation} 
\Se = \{ \U(i,j): \text{$j$ odd} \} 
\cup \{ \Phi(i,j,\tilde i, \tilde j)\}. 
\end{equation} 

\begin{Theorem} \sl 
The subspace of $S_3(S_d)$ generated by all the coefficients of 
all the elements in $\Se$ equals $(I_X)_3$. \emph{A fortiori}, a binary 
$d$-ic $F$ can be written as the $e$-th power of a quadratic form, 
iff all the elements of $\Se$ vanish on $F$. 
\end{Theorem} 
The second part of the theorem has a more classical flavour. It is an 
instance of the theme that, any property of an algebraic form which is 
invariant under a change of coordinates can be characterised by the 
vanishing (or non-vanishing) of a set of concomitants. 

\smallskip 

\demo 
Firstly we show that all elements in  $\Se$ vanish at a general 
point $F \in X$. By a change 
of variable, we may assume that $F = Q^e$, where $Q = x_0 \, x_1$. 
By formula (\ref{Lemmadixon1}), 
\[ (F,F)_{2i} = (-1)^i \, \cN^{\, \rm II}_{e,e,i} \; Q^{2e-2i}, \] 
hence $\E(i,j) = 0$ for $j$ odd by Proposition~\ref{prop.hyper}.
If $j= 2k$, then 
\[ (Q^{2e-2i},F)_{2k} = (-1)^k \, 
\cN^{\, \rm II}_{2e-2i,e,i} \; Q^{3e-2i-2k}, \] 
which implies $\E(i,j) = \mu_{i,j} \, Q^{3e-2i-j}$. 
Hence $\Phi(i,j,\tilde i, \tilde j)=0$ by definition. 

Let $J$ denote the subspace of $(I_X)_3$ generated by the 
coefficients of all the elements in $\Se$. So far we have shown that 
$J \subseteq (I_X)_3$. Let $S_m \subseteq (I_X)_3$ be an irreducible 
submodule, and $\Psi$ the corresponding covariant. We have to 
show that the coefficients of $\Psi$ are in $J$. Now it is a classical 
result that each covariant of binary forms is a linear combination of 
iterated transvectants of $F$ (see \cite[\S 86]{GrYo}). Since $\Psi$ is 
a cubic covariant, 
\begin{equation} 
\Psi = \sum\limits_{3d-4i-2j=m} q_{i,j} \, \E(i,j) 
\label{psi.sum} \end{equation}
for some $q_{i,j} \in \complex$. By hypothesis $\Psi$ vanishes on $X$. 
Hence if (\ref{psi.sum}) involves only one summand, then it must come from an 
odd $j$ and the claim follows. 
Alternately assume that $q_{i,j},q_{\, \tilde i,\tilde j} \neq 0$, then 
the covariant 
$\mu_{i,j} \, \Psi + q_{\, \tilde i,\tilde j} \,\Phi(i,j,\tilde i,\tilde j)$
involves at least one fewer summand and vanishes on $X$. 
Hence we are done by induction. \qed 

\medskip 

In general $\Se$ is not the smallest set which would make this theorem true. 
Given a particular value of $d$, it can be pared down substantially using 
properties of covariants specific to $d$. 
\begin{Example} \rm 
Assume $d=8$. Now 
\[ \begin{aligned}
S_3(S_8) & = S_{24} \oplus S_{20} \oplus S_{18} \oplus S_{16} \oplus 
S_{14} \oplus S_{12}^2 \oplus S_{10} \oplus S_8^2 \oplus S_4 \oplus S_0, \\
S_2(S_{12}) & = S_{24} \oplus S_{20} \oplus S_{16} \oplus S_{12} \oplus 
S_8 \oplus S_6 \oplus S_4 \oplus S_0, 
\end{aligned} \]
hence by Corollary~\ref{gr}, 
\[ (I_X)_3 = S_{18} \oplus S_{14} \oplus S_{12} \oplus S_{10} 
\oplus S_8 \oplus S_6. \] 
Hence it will suffice to choose a subset $\Te \subseteq \Se$, such that 
$\Te$ contains only one covariant each of orders $\{18,14,12,10,8,6\}$. 
Now, for instance, $\E(0,3) \in \Se$ is a covariant of order $18$ which can 
be chosen as an element of $\Te$. (Of course, $\E(1,1)$ would do as well. 
Observe that $S_3(S_8)$ contains only one copy of $S_{18}$, this implies that
$\E(0,3)$ and $\E(1,1)$ are constant multiples of each other 
for a generic $F$.) Similarly we select $\Phi(0,6,1,4)$ as 
the order $12$ covariant. Continuing in this way, we may let 
\[ \Te = \{\E(0,3),\E(0,5),\E(0,7),\Phi(0,6,1,4),\Phi(0,8,1,6),\E(3,3) \}, \] 
and then the previous theorem is true verbatim with $\Se$ replaced by $\Te$. 
Thus, a binary octavic $F$ is the fourth power of a quadratic form iff 
the covariants 
\[ \begin{array}{lll} 
(F^2,F)_3, & (F^2,F)_5, & 13 \, (F^2,F)_6 - 63 \, ((F,F)_2,F)_4, \\ 
(F^2,F)_7, & ((F,F)_6,F)_3, & 195 \, (F^2,F)_8 - 2744 \, ((F,F)_2,F)_6
\end{array} \]
are zero. We do not know if $\Te$ can be shortened any further while retaining this 
property. 

Throughout, we have tacitly assumed that none of the elements in 
$\Te$ vanishes identically. This can be checked by a simple 
direct calculation, e.g., by specializing $F$ to $x_0^5 \, x_1^3$. 
\end{Example} 
%%%%%%%%%%%%%%%%%%%%%%%%%%%%%%%%%%%%%%%%%%%%%%%%%%%%%%%%%
\section{A note on terminology and history}
In this paper we have adopted the term `Feynman diagrams' 
following the usage of theoretical physicists. 
However, the historical roots of this notion, especially 
in the context of invariant theory, substantially predate Feynman's work.

Feynman diagrams, as known to physicists, seem to have first 
appeared in print in the work of Dyson~\cite{Dyson}, who accredits them to the 
unpublished work of Richard Feynman. 
However, the idea of using discrete combinatorial structures 
(for instance tree graphs) to describe the outcome of repeated
applications of differential operators goes back to A.~Cayley~\cite{Cayley2}. 
Classically, such a diagrammatic approach was used in invariant theory 
by Sylvester~\cite{Sylvester}, Clifford~\cite{Clifford} and 
Kempe~\cite{Kempe}. It is remarkable that Clifford used what would
now be called Fermionic or Berezin integration in order to 
explain the translation from graphs to actual covariants.
The diagrams which we have used here directly mirror  
the classical symbolic notation: arrows correspond to bracket factors,
and each vertex corresponds to a symbolic letter, to be repeated
as many times as the degree of the vertex. 
The formalism used in Olver and Shakiban~\cite{OShakiban} 
is somewhat different due to a normal ordering procedure 
inspired by Gelfand~\cite{GelfandD}; it is also explained 
in~\cite[Ch.~6]{Olver}. A generally excellent account of the history 
of the diagrammatic notation in physics and group theory can be 
found in \cite[Chapter 4]{Cvitanovic}. Finally, the interesting  pedagogical 
work of computer graphics pioneer J.~F.~Blinn~\cite{Blinn},
(who was inspired by Stedman's work~\cite{Stedman}) deserves 
to be mentioned.

\medskip 

\noindent {{\sc Acknowledgements:} \smaller 
The first author would like to express his gratitude to
Professors David Brydges and Joel Feldman for the invitation to 
visit the University of British Columbia. He gratefully
acknowledges the support of the Centre National de la Recherche
Scientifique during the academic year 2002-2003.
The second author would like to thank Professor James Carrell for 
his invitation to visit UBC. 
We are indebted to the authors of the package Macaulay-2,
the digital libraries maintained by the Universities of 
Cornell, G{\"o}ttingen and Michigan, as well as J-Stor and 
Project Gutenberg.} 
%%%%%%%%%%%%%%%%%%%%%%%%%%%%%%%%%%%%%%%%%%%%%%%%%%%%

%%%%%%%%%%%%%%%%%%%%%%%%%%%%

\pagebreak 
\parbox{6cm}{\small 
{\sc Abdelmalek Abdesselam} \\ 
Department of Mathematics\\
University of British Columbia\\
1984 Mathematics Road \\
Vancouver, BC V6T 1Z2 \\ Canada. \\ 
{\tt abdessel@math.ubc.ca}} 
\hfill 
\parbox{5cm}{\small 
LAGA, Institut Galil\'ee \\ CNRS UMR 7539\\
Universit{\'e} Paris XIII\\
99 Avenue J.B. Cl{\'e}ment\\
F93430 Villetaneuse \\ France.} 

\vspace{1.5cm} 

\parbox{6cm}{\small 
{\sc Jaydeep Chipalkatti} \\ 
Department of Mathematics\\
University of Manitoba \\ 
433 Machray Hall \\ 
Winnipeg MB R3T 2N2 \\ Canada. \\ 
{\tt chipalka@cc.umanitoba.ca}}
\end{document}